\documentclass[twocolumn]{autart}    % Enable thiand disable the
\usepackage{graphicx}          % Include this line if your
\usepackage{hyperref}                               % you prefer.
\usepackage{subfigure}
\usepackage{mathtools,array}
\usepackage{mysymbol}
\usepackage{flushend}
\usepackage{har2nat,natbib}
\usepackage{wrapfig}
\DeclarePairedDelimiter\norm{\lVert}{\rVert}

\newcommand{\set}[1]{\left\{#1\right\}}

\newcommand{\s}{\bar{s}}
\newcommand{\x}[2]{{\mathbf{#1}}{(#2)}}
\newcommand{\xx}[3]{{\mathbf{#1}}_{#2}{(#3)}}
\newcommand{\supp}{\textup{supp}}
\newcommand{\rr}[2]{\mathbb{#1}^{#2}}
\usepackage{multirow}
\usepackage{amsmath}
\usepackage{amssymb}

\usepackage{algorithm}
\usepackage[noend]{algpseudocode}
\algrenewcommand\algorithmicforall{\textbf{foreach}}
\usepackage{etoolbox}

\errorcontextlines\maxdimen

\makeatletter
\newcommand*{\algrule}[1][\algorithmicindent]{\makebox[#1][l]{\hspace*{.5em}\vrule height .75\baselineskip depth .25\baselineskip}}%

\newcount\ALG@printindent@tempcnta
\def\ALG@printindent{%
    \ifnum \theALG@nested>0% is there anything to print
        \ifx\ALG@text\ALG@x@notext% is this an end group without any text?
            \addvspace{-3pt}% FUDGE for cases where no text is shown, to make the rules line up
        \else
            \unskip
            \ALG@printindent@tempcnta=1
            \loop
                \algrule[\csname ALG@ind@\the\ALG@printindent@tempcnta\endcsname]%
                \advance \ALG@printindent@tempcnta 1
            \ifnum \ALG@printindent@tempcnta<\numexpr\theALG@nested+1\relax% can't do <=, so add one to RHS and use < instead
            \repeat
        \fi
    \fi
    }%
\usepackage{etoolbox}
\patchcmd{\ALG@doentity}{\noindent\hskip\ALG@tlm}{\ALG@printindent}{}{\errmessage{failed to patch}}
\makeatother
\usepackage{booktabs}

\usepackage{mathtools}
\DeclarePairedDelimiter\ceil{\lceil}{\rceil}

\newcommand{\vz}{$v_l^0$ }
\newcommand{\vo}{$v_l^1$ }

\makeatletter
\newsavebox\myboxA
\newsavebox\myboxB
\newlength\mylenA

\newcommand*\xoverline[2][0.9]{%
    \sbox{\myboxA}{$\m@th#2$}%
    \setbox\myboxB\null% Phantom box
    \ht\myboxB=\ht\myboxA%
    \dp\myboxB=\dp\myboxA%
    \wd\myboxB=#1\wd\myboxA% Scale phantom
    \sbox\myboxB{$\m@th\overline{\copy\myboxB}$}%  Overlined phantom
    \setlength\mylenA{\the\wd\myboxA}%   calc width diff
    \addtolength\mylenA{-\the\wd\myboxB}%
    \ifdim\wd\myboxB<\wd\myboxA%
       \rlap{\hskip 0.5\mylenA\usebox\myboxB}{\usebox\myboxA}%
    \else
        \hskip -0.5\mylenA\rlap{\usebox\myboxA}{\hskip 0.5\mylenA\usebox\myboxB}%
    \fi}
\makeatother

\usepackage{xcolor}

\begin{document}

\begin{frontmatter}
%\runtitle{Insert a suggested running title}  % Running title for regular
                                              % papers but only if the title
                                              % is over 5 words. Running title
                                              % is not shown in output.

  \title{{An Optimal Graph-Search Method for  Secure State Estimation}\thanksref{footnoteinfo}} % Title, preferably not more % than 10 words.

\thanks[footnoteinfo]{{This work is supported in part by the ONR under agreements
    \#N00014-18-1-2374 and \#N00014-17-1-2504, the AFOSR under award \#FA9550-19-1-0169, as well as the NSF under grant CNS-1652544.} This paper was not presented at any IFAC meeting.}
\thanks[footnoteinfoauthor]{Corresponding author.}
\author[a]{Xusheng Luo \thanksref{footnoteinfoauthor}}\ead{xusheng.luo@duke.edu},    % Add the
\author[b]{Miroslav Pajic}\ead{miroslav.pajic@duke.edu},               % e-mail address
\author[a]{Michael M. Zavlanos}\ead{michael.zavlanos@duke.edu}  % (ead) as shown

\address[a]{Department of Mechanical Engineering and Materials Science, Duke University, Durham, NC 27708, U.S.A}
\address[b]{Department of Electrical and Computer Engineering, Duke University, Durham, NC 27708, U.S.A}

                             % keyword list or with the
                                          % help of the Automatica
                                          % keyword wiz
\begin{abstract}                          % Abstract of not more than 200 words.
  The growing complexity of modern Cyber-Physical Systems (CPS) and the frequent communication between their components make them vulnerable to malicious attacks. As a result, secure state estimation is a critical requirement for the control of these systems. Many existing secure state estimation methods suffer from combinatorial complexity which grows with the number of states and sensors in the system. This complexity can be mitigated using optimization-based methods that relax the original state estimation problem, although at the cost of optimality as these methods often identify attack-free sensors as attacked. {In this paper, we propose a new optimal graph-search algorithm  to correctly identify malicious attacks and to securely estimate the states even in large-scale CPS modeled as linear time-invariant systems.} The graph consists of layers, each one containing two nodes capturing a truth assignment of any given sensor, and directed edges connecting adjacent layers only. {Then, our algorithm searches the layers of this graph incrementally, favoring directions at higher layers with more attack-free assignments, while actively managing a repository of nodes to be expanded at later iterations. The proposed search bias and the ability to revisit nodes in the repository and self-correct,  allow our graph-search algorithm to reach the optimal assignment {faster} and tackle larger problems.} We show that our algorithm is complete and optimal provided that process and measurement noises do not dominate the attack signal. Moreover, we provide numerical simulations that demonstrate the ability of our algorithm to correctly identify attacked sensors and securely reconstruct the state. Our simulations show that our method outperforms existing algorithms both in terms of optimality and execution time.
\end{abstract}

\end{frontmatter}
\endNoHyper % Magic!
\section{Introduction}

Cyber-Physical Systems (CPS) are networked systems consisting of embedded physical components, such as sensors and actuators, and computational components, such as controllers. Recently, CPS have been successfully utilized in large-scale applications, including power network control, industrial manufacturing processes, and traffic control. However, the growing complexity of CPS and the frequent communication between  components make them vulnerable to malicious attacks. Such attacks often manipulate the state of the system by injecting faulty data through compromised sensors, leading to undesirable feedback control signals. Recently, cyber-attacks have been responsible for some incidents of safety-critical automobiles \cite{koscher2010experimental,auto,shoukry2013non} and UAVs \cite{javaid2017analysis}, and even worse, catastrophic losses of large-scale systems \cite{slay2007lessons,chen2011lessons}. Therefore, developing attack-resilient methods for secure state estimation in CPS has recently gained significant attention \cite{lee2008cyber,cardenas2008secure}.

In this paper, we consider CPS modeled as linear time-invariant systems, where a subset of sensors is subject to malicious attacks represented as attack vectors added to the measurements. Our goal is to detect the attacked sensors fast and use the attack-free sensors to accurately estimate the state. In this context, secure state estimation is closely related to robust control \cite{pasqualetti2011cyber,manandhar2014combating}, where the control design is subject to process and measurement noise, modeled as an unknown disturbance that is bounded or follows some probability distribution. Nevertheless, such assumptions on the noise restrict the general application of robust control methods for secure state estimation, since it is difficult to predict the  attack strategy. Similarly, fault tolerant control methods \cite{blanke2006diagnosis,teixeira2015secure} focus on internal faults with known failure modes and statistical properties, rather than  adversarial attacks. Secure state estimation under specific attack signals has been investigated in \citet{teixeira2010cyber,sundaram2010wireless,teixeira2012attack,miao2013stochastic,mo2014resilient,hendrickx2014efficient,mo2014detecting}.

{Compared to the literature discussed above, we do not impose any assumptions on the type of the attack signal.} We assume that the number of attacked sensors is {smaller than an upper bound}, which is necessary to ensure observability of the attacked system that is needed to reconstruct the state. Under this assumption, we propose a new optimal graph search-based algorithm to correctly identify malicious attacks {even in large-scale CPS} and securely estimate their state, {when the power of attack signals exceeds a certain threshold.} The graph consists of layers, each one containing two nodes capturing a truth assignment of any given sensor, and directed edges connecting adjacent layers only. Then, our algorithm searches the layers of this graph incrementally, favoring directions with more attack-free assignments and higher layer, while actively managing a repository of nodes whose expansion are intentionally delayed. The combination of search bias, intentionally delayed expansion and the ability to self-correct allow our graph-search algorithm to reach the optimal assignment fast and tackle larger problems. Assuming that process and measurement noise does not dominate the attack signal, we show that our algorithm is complete and optimal meaning that it will find a feasible attack assignment, if one exists, which does not incorrectly identify any  attack-free sensor as attacked. Finally, numerical simulations show that our method outperforms existing algorithms both in terms of optimality and execution time.

Most closely related to the work proposed here are the methods in \citet{fawzi2014secure,pajic2014robustness,pajic2015attack,pajic2017attack,chong2015observability,shoukry2017secure,mishra2017secure,shoukry2018smt}, which exploit the measurement within a finite-length time window to conduct state estimation. Specifically, \citet{fawzi2014secure} consider discrete-time LTI systems without noise and provides necessary and sufficient conditions under which the state of the system can be reconstructed when a subset of the sensors are under attack. The idea is to formulate the secure state estimation problem as an $\ell_0$ minimization problem that is computationally expensive, and then relax it into an $\ell_1/ \ell_r$ problem that can be efficiently solved using convex optimization, which mitigates the combinatorial complexity of the methods in \citet{pasqualetti2013attack,yong2015resilient,lee2015secure}. However, $\ell_0$ and $\ell_1/\ell_r$ optimization are not always equivalent, thus their relaxation can result in incorrect estimates. The work in~\citet{pajic2014robustness} extends this method to LTI systems where process and measurement noises are considered in the presence of malicious attacks, and formulates the $\ell_0$ optimization problem  as a mixed integer linear program (MILP). {However, solving MILPs is NP-hard, so this method can not be used for very large problems.}  Analytic bounds on the state-estimation error for the proposed $\ell_0$ state estimator and its convex $\ell_1$ relaxation in the presence of noise are derived in~\citet{pajic2015attack,pajic2017attack}, where it is shown that using relaxation results in inaccurate estimation. The work in \citet{chong2015observability} provides similar necessary and sufficient conditions for continuous-time LTI systems, and proposes two methods to estimate the state involving the observability Gramian and the Luenberger observer. However, both methods are computationally expensive. An alternative approach   based on Satisfiability Modulo Theory (SMT) is proposed in \citet{shoukry2017secure,mishra2017secure,shoukry2018smt}, that formulates the secure state estimation problem as a satisfiability problem subject to  Boolean constraints and convex constraints over real state variables. The proposed iterative algorithm
combines SMT solvers to obtain a possible attack assignment for the sensors with convex optimization methods to check whether this  assignment is valid given the dynamical system equations. %% If not, certificates are generated to impose additional constraints on the SMT solver, until a feasible solution is obtained.
Due to the formulation as a feasibility problem, the solution is not guaranteed to be optimal even in the absence of process and measurement noise. Compared to the literature discussed above, our graph-search method is provably optimal, meaning it identifies the true attack assignment and does not incorrectly identify attack-free sensors as attacked. Moreover, numerical experiments show that our method compares favorably to existing methods in terms of execution time. This is due to the proposed search bias that favors directions at higher layers with more attack-free assignments and the ability of our algorithm to self-correct by managing a repository of nodes that can be expanded at later iterations if needed.

The rest of the paper is organized as follows. Section~\ref{sec:prob} provides the problem formulation. In Section~\ref{sec:search}, we present the proposed graph-search algorithm for secure state estimation, and examine its completeness, optimality, and complexity in Section~\ref{sec:opt}. Finally, comparative numerical simulations are shown in Section~\ref{sec:sim}, while Section~\ref{sec:conclusion} concludes the paper.

\section{Problem Formulation}\label{sec:prob}

%% Given $p$ column vectors of the same dimension $\mathbf{y}_1, \ldots, \mathbf{y}_p \in \mathbb{R}^{T}, T \in \mathbb{N}^+$, we stack them together to constitute a block vector denoted by $\mathbf{Y} = [\mathbf{y}_1^T, \ldots, \mathbf{y}_p^T]^T \in \mathbb{R}^{pT}$. Similar to the support of a vector, the support of a block vector, denoted by $\supp(\mathbf{Y}) $, is defined as the set of indices of non-zero block vectors in $Y$, i.e., $\supp(\mathbf{Y}) = \{i \,|\, \mathbf{y}_i \not= \mathbf{0} \}\subseteq \{1, \ldots, p\}$. Similarly, we can define a block matrix $M$ made up of $T$ matrices $M_1, \ldots, M_{T} \in \mathbb{R}^{p\times n}$ as $M = [M_1^T, \ldots, M_p^T]^T \in \mathbb{R}^{pT\times n}$. Then, the support of a block matrix is the set of indices of non-zero matrices.

 %% For a given set $\Gamma\subseteq\{1,\ldots, p\}$, the symbols $\mathbf{Y}_{\Gamma}$ and $M_{\Gamma}$, respectively, represent the block vector and block matrix obtained after removing from $\mathbf{Y}$ and $M$ corresponding vectors and matrices indexed by $\Gamma$.

\subsection{Linear Dynamical Systems under Attack}
Consider the linear time-invariant dynamical system:
\begingroup\makeatletter\def\f@size{10}\check@mathfonts
\def\maketag@@@#1{\hbox{\m@th\normalsize\normalfont#1}}%
\begin{equation}
\begin{aligned}\label{equ:sys}
  \x{x}{t+1} &= A\x{x}{t} + B\x{u}{t} + \x{v}{t}, \\
  \x{y}{t} & = C\x{x}{t} + \x{e}{t} + \x{w}{t}.
\end{aligned}
\end{equation}
\endgroup
where $\x{x}{t}\in\rr{R}{n}, \x{u}{t}\in\rr{R}{m}$, and $ \x{y}{t}\in \rr{R}{p}$ denote the state vector, control vector, and measurement vector for $p$ sensors at time instant $t$, respectively; $A, B, C$ are system matrices with appropriate dimensions; $\x{v}{t}$ and $\x{w}{t}$ represent the process noise and measurement noise at time $t$; and $\x{e}{t}\in \rr{R}{p}$ is an attack vector so that if the $i$-th element $\xx{e}{i}{t}$ of $\x{e}{t}$ is non-zero then sensor $i$ is under attack, and is attack-free otherwise. We assume that the set of sensors that the attacker has access to does not change over time. Moreover, let $|\supp(\x{e}{t})|$ denote the number of attacked sensors, where $\supp(\mathbf{e}(t))\subseteq \{1,...,p \}$ denotes the support of the vector $\mathbf{e}(t)\in \rr{R}{p}$, that is the set of indices that correspond to non-zero elements in $\mathbf{e}(t)$, $|\cdot|$ denotes the cardinality of a set.
%% which are assumed to be upper bounded by constants $\xoverline{\mu}$ and $\xoverline{\psi}$, respectively, that is, $\norm{\x{\mu}{t}} \leq \xoverline{\mu}$ and $\norm{\x{w}{t}}\leq \xoverline{\psi}, \forall t \in \rr{N}{}$. Finally,\

In this paper, we do not consider the case where the actuators are under attack, thus, all the control inputs are known and we can subtract their effect from the dynamical equations due to linearity. Therefore, for simplicity, we set the matrix $B$ to be zero. Given $T$ measurements $\x{y}{t-T+1}, \ldots, \x{y}{t}$ that are subject to attack vectors $\x{e}{t-T+1}, \ldots, \x{e}{t}$, we can express them as a function of the  states $\x{x}{t-T+1}$ as:
\begingroup\makeatletter\def\f@size{10}\check@mathfonts
\def\maketag@@@#1{\hbox{\m@th\normalsize\normalfont#1}}%
\begin{align}\label{equ:meas}
  \xx{Y}{t,T}{t} & = \ccalO_T \,\x{x}{t-T+1} + \xx{e}{t,T}{t} + \xx{w}{t,T}{t},
\end{align}
where
\begin{align}
  \xx{Y}{t,T}{t} & = \left[\x{y^\intercal}{t-T+1}, \, \x{y^\intercal}{t-T+2}, \, \cdots, \,  \x{y^\intercal}{t}\right]^\intercal, \nonumber \\
  \ccalO_T \;\, & = \left[ C^\intercal, \, A^\intercal C^\intercal, \, \cdots, \,  (A^\intercal)^{T-1} C^\intercal \right]^\intercal, \nonumber \\
  \xx{e}{t,T}{t} & = \left[  \x{e^\intercal}{t-T+1}, \, \x{e^\intercal}{t-T+2}, \,  \cdots, \, \x{e^\intercal}{t} \right]^\intercal, \nonumber \\
  \xx{w}{t,T}{t} & = \left[
    \begin{array}{c}
      \x{w}{t-T+1}\\
      C \x{v}{t-T+1} + \x{w}{t-T+2}\\
       C\sum_{i=1}^{2}A^{2-i} \x{v}{t-T+i} + \x{w}{t-T+3}\\  %C\left[A^2 \x{v}{t-T+1}+ A \x{v}{t-T+2}\right] + \x{w}{t-T+3} \\
      \ldots\\
      C\sum_{i=1}^{T-1}A^{T-1-i} \x{v}{t-T+i} + \x{w}{t}
    \end{array}
    \right]. \nonumber
\end{align}
\endgroup
%% \begin{align*}
%%   \xx{Y}{t,T}{t} = \left[\x{y}{t-T+1}, \x{y}{t-T+2}, \cdots, \x{y}{t}\right],\,
%%   \xx{Y}{t,T}{t} = \left[
%%     \begin{array}{@{\mkern-0.1mu} c @{\mkern-0.1mu}}
%%       \x{y}{t-T+1}\\
%%       \x{y}{t-T+2}\\
%%       \vdots\\
%%       \x{y}{t}
%%     \end{array}
%%     \right],\,
%%   \ccalO_T = \left[
%%     \begin{array}{@{\mkern-0.1mu} c @{\mkern-0.1mu}}
%%     C\\
%%     CA\\
%%       \vdots\\
%%       CA^{T-1}
%%     \end{array}
%%     \right],
%% \end{align*}
%% \begin{align*}
%%    \xx{e}{t,T}{t} = \left[
%%     \begin{array}{@{\mkern-0.1mu} c @{\mkern-0.1mu}}
%%       \x{e}{t-T+1}\\
%%       \x{e}{t-T+2}\\
%%       \vdots\\
%%       \x{e}{t}
%%     \end{array}
%%     \right],
%% \end{align*}
%% and
%% \begin{align*}
%%   \xx{w}{t,T}{t} = \left[
%%     \begin{array}{c}
%%       \x{w}{t-T+1}\\
%%       C \x{v}{t-T+1} + \x{w}{t-T+2}\\
%%        C\sum_{i=1}^{2}A^{2-i} \x{v}{t-T+i} + \x{w}{t-T+3}\\  %C\left[A^2 \x{v}{t-T+1}+ A \x{v}{t-T+2}\right] + \x{w}{t-T+3} \\
%%       \vdots\\
%%       C\sum_{i=1}^{T-1}A^{T-1-i} \x{v}{t-T+i} + \x{w}{t}
%%     \end{array}
%%     \right]
%% \end{align*}
The term $\xx{e}{t,T}{t}$ denotes the attack vector, and with a slight abuse of notation, $\xx{w}{t,T}{t}$ represents noise vectors, including both the process and measurement noise. If the system is noiseless and attack-free, equation~\eqref{equ:meas} becomes
\begingroup\makeatletter\def\f@size{10}\check@mathfonts
\def\maketag@@@#1{\hbox{\m@th\normalsize\normalfont#1}}%
\begin{align}\label{equ:noiseless}
\xx{Y}{t,T}{t} = \ccalO_T \,\x{x}{t-T+1}.
\end{align}
\endgroup
As is shown in \citet{fawzi2014secure}, for noiseless systems that are under attack, we can reconstruct the state from $T$ measurements when $s$ sensors are under attack if and only if, $ \forall\, \mathbf{x}\in\rr{R}{n}\setminus\set{\mathbf{0}}$,
$| \supp(C\mathbf{x}) \cup  \supp(CA\mathbf{x}) \cup\ldots\cup \,  \supp(CA^{T-1} \mathbf{x})  |  > 2s.$
Therefore, there is an upper limit on the number of attacked sensors, denoted by $\bar{s}$, beyond which states can not be correctly estimated. This limit depends on the system matrices $A$ and $C$. In this paper, we assume the number of attacked sensors is less than or equal to~$\s$, i.e., $|\supp(\x{e}{t})| \leq \s$, and $\s$ is assumed to be known a priori, a common assumption used in relevant work. {Furthermore, we assume the system in~\eqref{equ:sys} is 2$\s$-sparse observable, which means the system is still observable after any $2\s$ sensors are removed. As shown in \citet{fawzi2014secure,shoukry2017secure}, it is impossible to correctly estimate the state if  $\ceil*{p/2}$ or more sensors are attacked. Thus, $\s \leq \ceil*{p/2}-1$. Besides these assumptions, we do not impose additional constraints on the attack vector, which can be arbitrary and unbounded.}

Moreover, let $\mathcal{I}\subseteq \{1,\dots ,p\}$ be a subset of sensors and define by $\xx{y}{}{t}|_{\ccalI}$ the vector composed of elements of $\mathbf{y}(t)$ indexed by the set $\ccalI$. Then, we can define the set of measurements corresponding to sensors in the set $\ccalI$ as
\begingroup\makeatletter\def\f@size{9}\check@mathfonts
\def\maketag@@@#1{\hbox{\m@th\normalsize\normalfont#1}}%
\begin{align*}
   \xx{Y}{t,T}{t}|_{\ccalI} = \left[\x{y^\intercal}{t-T+1}|_{\ccalI},
      \x{y^\intercal}{t-T+2}|_{\ccalI},   \cdots,  \x{y^\intercal}{t}|_{\ccalI} \right]^\intercal.
\end{align*}
\endgroup
Considering only measurements from sensors indexed by the set $\ccalI$,  equation \eqref{equ:meas} can be rewritten as
\begin{align*}
 \xx{Y}{t,T}{t}|_{\ccalI} = \ccalO_T|_{\ccalI} \,\x{x}{t-T+1} + \xx{e}{t,T}{t}|_{\ccalI} + \xx{w}{t,T}{t}|_{\ccalI}.
\end{align*}
For notational simplicity, we use  $\mathbf{Y}_{\ccalI}, \ccalO_{\ccalI}, \mathbf{e}_{\ccalI}, \mathbf{w}_{\ccalI}$ to denote  $ \xx{Y}{t,T}{t}|_{\ccalI}$, $\ccalO_T|_{\ccalI}, \xx{e}{t,T}{t}|_{\ccalI}, \xx{w}{t,T}{t}|_{\ccalI}$ and $\mathbf{x}$ to denote $\mathbf{x}(t-T+1)$ when they are clear from the context. Moreover, when $\ccalI = \set{i}$, for $i\in\set{1,\ldots,p}$, is a singleton, we use the notations $\mathbf{Y}_i, \ccalO_i, \mathbf{e}_i, \mathbf{w}_i$. Finally, we assume that the noise term $\xx{w}{t,T}{t}$ is upper bounded. This is a reasonable assumption since, otherwise, it is impossible to estimate the state. Specifically, we assume that $\norm{\mathbf{w}_i} \leq \xoverline{w}_i$, for $\forall t\in\rr{N}{}, T\in \set{1,\ldots,n}, i \in \set{1,\ldots,p}$, where $\norm{\cdot}$ is $\ell_2$-norm of a vector. Similarly, $\norm{\mathbf{w}_{\ccalI}}^2 \leq  \xoverline{w}_{\ccalI}^2  = \sum_{i\in\ccalI} \xoverline{w}_i^2$ and $\xoverline{w}^2 = \sum_{i=1}^p \xoverline{w}_i^2$.

\subsection{Secure State Estimation}
Let $\mathbf{b} = (b_1, \ldots, b_p)$ be a vector of binary variables such that $b_i = 1$ if sensor $i$ is under attack and $b_i=0$, otherwise. Let $\mathbf{x}_0$ and $ \mathbf{b}_0$ denote the true state and the true attack assignment at time $t-T+1$, respectively. Assuming the $|\supp{(\mathbf{b}_0)}|\leq \bar{s}$, our goal is to find the attack assignment $\mathbf{b}^*$ that satisfies $\mathbf{b}^* = \mathbf{b}_0$, and use the attack-free sensors in $\mathbf{b}^*$ to reconstruct the state. Specifically, we formulate the following problem.

\begin{prob}\label{pb:min}
   Consider the linear dynamical system in \eqref{equ:sys} that is under attack. Determine the optimal state vector and attack assignment $(\mathbf{x}^*,\mathbf{b}^*) \in \mathbb{R}^n \times \mathbb{B}^p$ that solve the optimization problem
 \begin{subequations}\label{equ:min}
\begin{equation}
 \hspace{-8em} \min_{(\mathbf{x},\mathbf{b}) \in \mathbb{R}^n \times \mathbb{B}^p} \quad |\supp(\mathbf{b})|
\tag{\ref{equ:min}}
\end{equation}
\vspace{-8mm}
\begin{align}
 s.t. & \quad\quad \| \mathbf{Y}_{\ccalI}- \ccalO_{\ccalI} \mathbf{x}\|_2 \leq \xoverline{w}_{\ccalI} + \sqrt{\epsilon}, \label{equ:mina}\\
& \quad\quad |\supp(\mathbf{b})| \leq \s, \label{equ:minb}
\end{align}
\end{subequations}
   where $\s$ is the maximum allowable number of attacked sensors in order to reconstruct states, $\ccalI = \{1,\ldots,p\}\setminus\supp(\mathbf{b})$, and $\epsilon$ is the acceptable numerical accuracy in the solution specified by the user, which serves as the stopping criterion of numerical iterations.
\end{prob}
As shown in \citet{shoukry2017secure}, if the  noise and solution accuracy are zero, i.e., if $\xoverline{w}_i=0$ and $\epsilon=0$, and if the system is $2\s$-sparse observable, then any assignment $\mathbf{b}$ with $\xoverline{\supp{(\mathbf{b})}} \subseteq \xoverline{\supp{(\mathbf{b}_0)}}$ and $|\supp{(\mathbf{b})}| \leq \s$ is a feasible assignment, where $\xoverline{\supp(\mathbf{b})}$ is the complement of the set $\supp(\mathbf{b})$, i.e., $\xoverline{\supp{(\mathbf{b})}} = \set{1,\ldots,n}\setminus\supp(\mathbf{b})$. In  words, a feasible solution to Problem~\ref{pb:min} correctly identifies all attacked sensors, but it can also incorrectly treat attack-free sensors as attacked. As shown in~\citet{shoukry2016event}, if the noise and solution accuracy are zero, then the solution to Problem~\ref{pb:min} correctly identifies the true attack assignment, i.e., it satisfies $\mathbf{b}^*=\mathbf{b}_0$. However, in the presence of noise and non-zero solution accuracy, a feasible solution to Problem~\ref{pb:min} can incorrectly identify attacked sensors as attack-free if the attack signal is undetectable meaning that its effect is  hidden by the noise and solution accuracy.\footnote{Typically, noise is the main factor that can hide the effect of an attack signal.} Similarly, {if the noise is relatively large or the solution accuracy is low, then,} some attack-free sensors may be incorrectly treated as being under attack. Nevertheless, if the attack vector is strong enough so that it can not be hidden by noise and solution accuracy, then it is reasonable to expect that the solution to Problem~\ref{pb:min} coincides with the true attack. {This discussion on attack detection in the presence of noise is also supported by the theoretical analysis in~\citet{pajic2017attack}. The following proposition quantifies this discussion in the current problem formulation.}

However, before we show this result, we provide some definitions. Let $\mathbf{x}'$ be any reachable state of system~\eqref{equ:sys}. Then, using the noiseless and attack-free model~\eqref{equ:noiseless}, we can get $T$ measurements $\mathbf{Y'} = \ccalO \mathbf{x}'$. Let $\tilde{\mathbf{x}}$ denote the solution of the problem $\min_{\mathbf{x}\in \rr{R}{n}} \norm{\mathbf{Y}'  - \ccalO \mathbf{x}}$ for the given $\mathbf{x}'$. Then, we can define the solution accuracy $\epsilon^*$ as
\begin{align}\label{equ:ep}
\epsilon^* = \inf\{ \epsilon\, |\, \norm{ \mathbf{Y}'  - \ccalO \tilde{\mathbf{x}}} \leq \sqrt{\epsilon}, \, \forall\, \mathbf{x}'\}.
\end{align}
In other words, $\epsilon^*$ is a uniform lower bound on $\epsilon$ so that for any $\mathbf{x}'$ the solution $\tilde{\mathbf{x}}$ of the problem $\min_{\mathbf{x}\in \rr{R}{n}} \norm{\mathbf{Y}'  - \ccalO \mathbf{x}}$ satisfies $\norm{\mathbf{Y}'  - \ccalO \tilde{\mathbf{x}}} \leq \sqrt{\epsilon^*}$. If $\epsilon < \epsilon^*$, then this constraint can not be satisfied. Now, consider a set $\mathcal{I}$ containing only attack-free sensors. Then, the solution $\bar{\mathbf{x}}$ to the problem  $\min_{\mathbf{x}\in \rr{R}{n}} \norm{\mathbf{Y}'_{\ccalI}  - \ccalO_{\ccalI} \mathbf{x}}$ also satisfies $\norm{\mathbf{Y}'_{\ccalI}  - \ccalO_{\ccalI} \bar{\mathbf{x}}} \leq \sqrt{\epsilon^*}$. The reason is that for $\tilde{\mathbf{x}}$, we have $\norm{\mathbf{Y}'_{\ccalI}  - \ccalO_{\ccalI} \tilde{\mathbf{x}}} \leq \norm{\mathbf{Y}'  - \ccalO \tilde{\mathbf{x}}}$, and $\norm{\mathbf{Y}'  - \ccalO \tilde{\mathbf{x}}} \leq \sqrt{\epsilon^*}$ by definition of $\epsilon^*$. Therefore,  $\norm{\mathbf{Y}'_{\ccalI}  - \ccalO_{\ccalI} \tilde{\mathbf{x}}} \leq \sqrt{\epsilon^*}$. Moreover, since $\bar{\mathbf{x}}$ is the minimizer of $\norm{\mathbf{Y}'_{\ccalI}-\ccalO_{\ccalI} \mathbf{x}}$, we have  $ \norm{\mathbf{Y}'_{\ccalI}-\ccalO_{\ccalI} \mathbf{\bar{x}}} \leq  \norm{ \mathbf{Y}'_{\ccalI}-\ccalO_{\ccalI}\tilde{\mathbf{x}}}$. We conclude that  $\norm{\mathbf{Y}'_{\ccalI}-\ccalO_{\ccalI} \bar{\mathbf{x}}} \leq \sqrt{\epsilon^*}$.

\begin{prop}\label{thm:attack}
Let the linear system in~\eqref{equ:sys} be $2\s$-sparse observable and let the number of attacked sensors be less than or equal to $\s$. Moreover, let $\epsilon=\epsilon^*$, i.e., the lower bound defined in~\eqref{equ:ep}. If the attack signal satisfies
\begin{align}\label{equ:attack}
  \norm{\mathbf{e}_{i}} > \left( \frac{2}{\sqrt{1-\Delta_s}}\right) \xoverline{w} + \frac{\sqrt{\epsilon}}{\sqrt{1-\Delta_s}},
\end{align}
where
  \begingroup\makeatletter\def\f@size{9}\check@mathfonts
\def\maketag@@@#1{\hbox{\m@th\normalsize\normalfont#1}}%
\begin{displaymath}
  \Delta_s  = \max_{\substack{\Gamma\subset\ccalI \subset \set{1,\ldots,p} \\ |\Gamma|\leq\s, |\ccalI|\geq p - \s}} \lambda_{\max} \left\{ \left( \sum_{i\in\Gamma} \ccalO_i^T \ccalO_i\right) \left( \sum_{i\in\ccalI} \ccalO_i^T \ccalO_i\right)^{-1} \right\},
\end{displaymath}
\endgroup
and $\lambda_{\max}(\cdot)$ is the maximal eigenvalue of a matrix, then $\mathbf{b}_0$ is the unique optimal solution of Problem~\ref{pb:min}.
\end{prop}

\begin{pf}
  First,  we show that if there is an attacked sensor in the index set $\ccalI$ with $|\ccalI| \ge p - \s$, when~\eqref{equ:attack} is satisfied, $\min_{\mathbf{x}\in \rr{R}{n}}\norm{\mathbf{Y}_{\ccalI} - \ccalO_{\ccalI} \mathbf{x}} \leq \xoverline{w}_{\ccalI} + \sqrt{\epsilon}$  does not hold anymore. It is shown in Theorem IV.3 in \citet{shoukry2017secure} that if there is an attacked sensor in the index set $\ccalI$ with $|\ccalI| \ge p - \s$, for which the attack signal satisfies~\eqref{equ:attack}, then the following two inequalities hold,
  \begingroup\makeatletter\def\f@size{9}\check@mathfonts
\def\maketag@@@#1{\hbox{\m@th\normalsize\normalfont#1}}%
  \begin{align}
    &  \norm{\mathbf{Y}_{\ccalI} - \ccalO_{\ccalI} \mathbf{x}}^2 \geq  {\left(\norm{(I -\ccalO_{\ccalI}\ccalO_{\ccalI}^+) \mathbf{e}_{\ccalI}} - \norm{(I -\ccalO_{\ccalI}\ccalO_{\ccalI}^+) \mathbf{w}_{\ccalI}}\right)}^2,  \label{eq:13} \\
   & \norm{(I -\ccalO_{\ccalI}\ccalO_{\ccalI}^+) \mathbf{e}_{\ccalI}} - \norm{(I -\ccalO_{\ccalI}\ccalO_{\ccalI}^+) \mathbf{w}_{\ccalI}} > \xoverline{w}_{\ccalI} + \sqrt{\epsilon}, \label{eq:14}
  \end{align}
  \endgroup
  %% \begin{align}\label{eq:14}
  %%   \norm{(I -\ccalO_{\ccalI}\ccalO_{\ccalI}^+) \mathbf{e}_{\ccalI}} - \norm{(I -\ccalO_{\ccalI}\ccalO_{\ccalI}^+) \mathbf{w}_{\ccalI}} & > \xoverline{w}_{\ccalI} + \sqrt{\epsilon},
  %% \end{align}
{where $\ccalO_{\ccalI}^+ = (\ccalO_{\ccalI}^T \ccalO_{\ccalI})^{-1} \ccalO_{\ccalI}^T $, $I$ is the identity matrix, and~\eqref{eq:13} corresponds to (13) and~\eqref{eq:14} is the second to last inequality in (14) in~\citet{shoukry2017secure}. It's worth noting that the right-hand side of constraint~\eqref{equ:mina} takes the form  $\xoverline{w}_{\ccalI} + \epsilon$ in equation $(8)$ in~\citet{shoukry2017secure}. {However, this reformulation does not invalidate the theoretical results in~\citet{shoukry2017secure}, since the derivation of inequalities~\eqref{eq:13} and~\eqref{eq:14} is irrelevant to any form of constraint~\eqref{equ:mina}.} Thus, combining~\eqref{eq:13} and~\eqref{eq:14}, we have that $\min_{\mathbf{x}\in \rr{R}{n}}\norm{\mathbf{Y}_{\ccalI} - \ccalO_{\ccalI} \mathbf{x}} \leq \xoverline{w}_{\ccalI} + \sqrt{\epsilon}$ does not hold anymore.} Therefore, the set $\ccalI$ should only include attack-free sensors. This implies that when inequality~\eqref{equ:attack} is satisfied, any feasible solution $\mathbf{b}$ of Problem~\ref{pb:min} needs to correctly identify all attacked sensors, i.e., $\supp{(\mathbf{b}_0)} \subseteq \supp{(\mathbf{b})}$, and therefore $|\supp{(\mathbf{b}_0)}| \leq |\supp{(\mathbf{b})}|$.

  In what follows, we show that the converse is also true, i.e., that any assignment $\mathbf{b}$ that satisfies $\supp(\mathbf{b}_0)\subseteq \supp(\mathbf{b})$ and $|\supp(\mathbf{b})| \leq \bar{s}$ is a feasible solution to Problem~\ref{pb:min}, provided that~\eqref{equ:attack} is satisfied. Then, we can conclude that since $\mathbf{b} = \mathbf{b}_0$ satisfies $\supp(\mathbf{b}_0)\subseteq \supp(\mathbf{b})$ with equality and $|\supp(\mathbf{b})| \leq \bar{s}$, $\mathbf{b}_0$ is a feasible and in fact the optimal solution to Problem~\ref{pb:min}. Moreover, since $\mathbf{b}_0$ is a binary vector, there is no other assignment $\mathbf{b}$ that satisfies $\supp(\mathbf{b}_0)\subseteq \supp(\mathbf{b})$ with equality and $|\supp(\mathbf{b})| \leq \bar{s}$. Therefore, $\mathbf{b}_0$ is the unique optimal solution to Problem~\ref{pb:min}. Note that if $\epsilon < \epsilon^*$, then it is possible that an assignment $\mathbf{b}$ satisfies $\supp(\mathbf{b}_0)\subseteq \supp(\mathbf{b})$ and $|\supp(\mathbf{b})| \leq \bar{s}$ but does not satisfy the constraint~\eqref{equ:mina}. Therefore, constraint~\eqref{equ:attack} is a necessary condition.

  To show that any assignment $\mathbf{b}$ that satisfies $\supp(\mathbf{b}_0)\subseteq \supp(\mathbf{b})$ and $|\supp(\mathbf{b})| \leq \bar{s}$ is a feasible solution to Problem~\ref{pb:min}, note first that if all sensors in the set $\ccalI = \xoverline{\supp(\mathbf{b})}$ are attack-free, we have that $\mathbf{Y}_{\ccalI} = \ccalO_{\ccalI} \mathbf{x}_0 + \mathbf{w}_{\ccalI}$ for the true $\mathbf{x}_0$. Moreover, using the true state $\mathbf{x}_0$, we can generate measurements $\mathbf{Y}'_{\ccalI} = \ccalO_{\ccalI} \mathbf{x}_0 $ according to the noiseless and attack-free model in~\eqref{equ:noiseless}. Combining those two equations we get $\mathbf{Y}_{\ccalI} = \mathbf{Y}'_{\ccalI} + \mathbf{w}_{\ccalI}$. By definition of $\epsilon^*$,
  \begingroup\makeatletter\def\f@size{9}\check@mathfonts
\def\maketag@@@#1{\hbox{\m@th\normalsize\normalfont#1}}%
  \begin{align}\label{equ:def}
  \min_{\mathbf{x} \in \rr{R}{n}} \norm{\mathbf{Y}^{}_{\ccalI} - \ccalO_{\ccalI} \mathbf{x} - \mathbf{w}_{\ccalI}}= \min_{\mathbf{x} \in \rr{R}{n}} \norm{\mathbf{Y}'_{\ccalI} - \ccalO_{\ccalI} \mathbf{x}} \leq \sqrt{\epsilon^*},
  \end{align}
  \endgroup
  where $\mathbf{Y}_{\ccalI} = \mathbf{Y}'_{\ccalI} + \mathbf{w}_{\ccalI}$. Since by the triangle inequality $ \norm{\mathbf{Y}_{\ccalI} - \ccalO_{\ccalI} \mathbf{x} } \leq \norm{\mathbf{Y}_{\ccalI} - \ccalO_{\ccalI} \mathbf{x} - \mathbf{w}_{\ccalI}} +\norm{\mathbf{w}_{\ccalI}}$, and $\mathbf{w}_{\ccalI}$ is irrelevant to the minimization over $\mathbf{x}$, we have
  \begin{align}\label{equ:tri}
     \min_{\mathbf{x} \in \rr{R}{n}} \norm{\mathbf{Y}_{\ccalI} - \ccalO_{\ccalI} \mathbf{x} }    \leq  \min_{\mathbf{x} \in \rr{R}{n}} \norm{\mathbf{Y}_{\ccalI} - \ccalO_{\ccalI} \mathbf{x} - \mathbf{w}_{\ccalI}} +  \norm{\mathbf{w}_{\ccalI}}.
  \end{align}
  Finally, combining~\eqref{equ:def} and~\eqref{equ:tri}, we get
$\min_{\mathbf{x} \in \rr{R}{n}} \norm{\mathbf{Y}_{\ccalI} - \ccalO_{\ccalI} \mathbf{x} }  \leq \norm{\mathbf{w}_{\ccalI}}  + \sqrt{\epsilon^*}.  $ %% \leq  \xoverline{w} + \sqrt{\epsilon^*}.
  Thus, providing~\eqref{equ:attack} holds, if the assignment $\mathbf{b}$ satisfies $\supp{(\mathbf{b}_0)} \subseteq \supp(\mathbf{b})$, then $\ccalI = \xoverline{\supp(\mathbf{b})}$ satisfies constraint~\eqref{equ:mina}, completing the proof.
\end{pf}

\section{Graph Search-based Secure State Estimation}\label{sec:search}
\begin{figure}[t]
  \centering
  \includegraphics[width=0.5\linewidth]{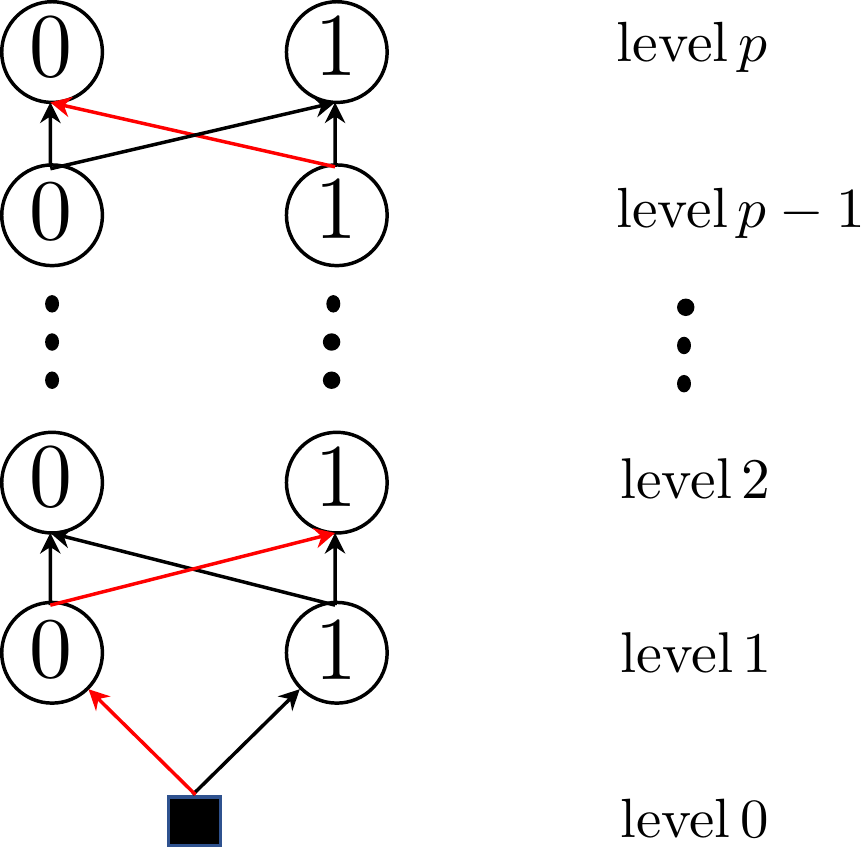}
  \caption{The red path stands for one possible full assignment $ \mathbf{b}= 0,1,\ldots,1,0$.}
  \label{fig:graph}
\end{figure}

In this section, we propose a graph search algorithm that incrementally assigns a truth value to each binary variable in $\mathbf{b}$. Specifically, our algorithm searches for the true attack assignment on a directed graph with $p+1$ levels and 2 nodes on each level; see Fig.~\ref{fig:graph}. The graph is initialized with an artificial root to provide a unique starting point.  This root corresponds to level 0. Each level except level 0 captures the truth assignment of one sensor, so that the nodes with values 1 and 0 at this level indicate whether this sensor is under attack or not, respectively. The edges in this graph connect nodes in adjacent levels only, and all edges point towards a higher level. A path starting from level 1 and ending at level $p$ corresponds to a possible attack assignment $\mathbf{b}$. Throughout the rest of this paper, we refer to the partial and full assignment when  part of or all of $p$ sensors are assigned, respectively.

{The key idea of the proposed algorithm is to prioritize search along paths in the graph that contain more attack-free sensors. Search along these paths returns sensor assignments with more 0 entries that minimize the objective in~\eqref{equ:min}. To decide whether a partial/full assignment is valid, our algorithm also checks the feasibility of the system of linear equations $\mathbf{Y}_\ccalI = \ccalO_{\ccalI} \mathbf{x}$, as per the constraints in~\eqref{equ:mina}.\footnote{Here, $\ccalI$ only includes sensors that have been assigned and assigned as attack-free.} During the early stages of the search, this system of equations may be underdetermined or square (when $T\cdot |\ccalI| \le n$), according to the observability matrix $\ccalO_\ccalI \in \mathbb{R}^{T|\ccalI|\times n}$. Therefore, it is possible that attacked sensors are incorrectly treated as attack-free. However, as the search progresses, more sensors are assigned 0 values, which reduces the dimension of the kernel space of the system of linear equations $\mathbf{Y}_\ccalI = \ccalO_{\ccalI} \mathbf{x}$, making the search less tolerant to such mistakes. To see this, assume that the observability matrix $\ccalO_{\ccalI}$ has full rank $n$ and the corresponding partial assignment is correct. Then, all subsequent sensors will also be correctly assigned. This is because, if attacked sensors are incorrectly identified as attack-free, then the system of linear equations $\mathbf{Y}_\ccalI = \ccalO_{\ccalI} \mathbf{x}$ will become inconsistent. On the other hand, if the observability matrix has full rank $n$ and the corresponding partial assignment is wrong, then the search along this path will terminate immediately as soon as the next sensor assigned 0 is added to the system of equations $\mathbf{Y}_\ccalI = \ccalO_{\ccalI} \mathbf{x}$, making it inconsistent. Together with prioritizing search along paths that contain more attack-free sensors to minimize the objective in~\eqref{equ:min}, the proposed algorithm also actively manages a repository of nodes whose exploration is deliberately postponed until there is need to correct a search that early on has made an incorrect sensor assignment. The combination of search bias and the ability to self-correct allow our graph-search algorithm to identify the true attack assignment fast.}   %% And also, a subpath with length $l$ gives us a partial assignment of the first $l$ sensors, thereby, we say subpath and partial assignment interchangeably, same as path and full assignment.

{Before presenting our algorithm, we introduce the information that a node needs to store and the resulting order between any two nodes. We first associate with every node a 5-tuple $(level, value, parent, \ccalI$, $residual)$, where:
(i) $level=l$ indicates that this node corresponds to the $l$-th sensor except when $l=0$ which corresponds to the artificial root;
(ii) $value=1$ means that the $l$-th sensor is under attack and $value=0$ means that the $l$-th sensor is attack-free;
(iii) $parent$ denotes the parent node of the $l$-th sensor obtained from level $l-1$;
(iv) $\ccalI$ is the set of levels with $value$ 0 (attack-free sensors) from level 1 to level $l$. The set $\ccalI$ and the value {\em level} together provide complete knowledge of the truth assignments of the first $l$ sensors, and they also contain information about the objective~\eqref{equ:min} and the satisfaction of constraint~\eqref{equ:minb}}; (v) $residual = 0$ if the inequality $\min_{\mathbf{x}\in\rr{R}{n}}\norm{\mathbf{Y}_{\ccalI} - \ccalO_{\ccalI}\mathbf{x}}\leq \xoverline{w}_{\ccalI} + \sqrt{\epsilon}$ holds, otherwise $residual=1$. {We can tell from $residual$ whether the current attack assignment of the first $l$ sensors satisfies constraint~\eqref{equ:mina}.}
Next, with the help of the 5-tuple, we define an ordering between two nodes $v$ and $v'$ in the graph as follows. We say $v = v'$ if $v.value = v'.value$ and $v.level = v'.level$. That is, two nodes are treated as equivalent if they correspond to the same sensor and have identical Boolean values, as shown in Fig.~\ref{fig:graph}. {Moreover, we use the lexicographical order to compare any two nodes.} Specifically, we define $v > v'$ if $v.level-|v.\ccalI| < v'.level - |v'.\ccalI|$ or $(v.level-|v.\ccalI| = v'.level - |v'.\ccalI|) \wedge (v.level > v'.level)$. In words, ordering of nodes is first determined by the number of attacked sensors in the paths leading to those nodes, so that fewer attacked sensors in these paths correspond to nodes with higher priority. This is because the goal of the algorithm is to find the true attack assignment and avoid incorrect assignments. On the other hand, if the paths leading to two nodes contain the same number of attacked sensors, ordering is determined by the levels of these nodes; higher level corresponds to higher priority in the ordering. This is because expanding a node at a higher level can advance the algorithm faster to terminate.
\begin{algorithm}[t]
  \caption{Graph-Search for Secure State Estimation}\label{alg:sse}
  %\vspace{0.5cm}
%\hspace*{\algorithmicindent} \textbf{Input:} measurements $\mathbf{Y}$, observability matrix $\ccalO$, \hspace*{18mm} number of sensors $p$,  noise bound $\xoverline{w}_i$, \hspace*{18.8mm} $ i \in \{1, \ldots, p\}$, solution accuracy $\epsilon$,\\
\hspace*{\algorithmicindent} \textbf{Input:} measurements $\mathbf{Y}$, observability matrix $\ccalO$,  number of sensors $p$,  noise bound $\xoverline{w}_i$,  $ i \in \{1, \ldots, p\}$, solution accuracy $\epsilon$, maximum allowable number of attacked sensors $\s$\\
\hspace*{\algorithmicindent} \textbf{Output:} estimated state $\mathbf{x}$, indices of attacked sensor $\ccalI$
\begin{algorithmic}[1]
  %% \KwIn {measurements $\mathbf{Y}$, observability matrix $\ccalO$, number of sensors $p$,  noise bound $\xoverline{w}_i, i \in \{1, \ldots, p\}$, solution accuracy $\epsilon$}
      %% \KwOut {estimated state $\mathbf{x}$, indices of attacked sensor $\ccalI$}
  \State$node.level \gets 0 ;\, node.value \gets1 ;\,  node.parent \gets  None;\,  node.\ccalI \gets \emptyset ;\,  node.residual = 0$\label{alg:line1}
  \State {\em frontier} $\gets\{node \}$; {\em explored} $\gets \emptyset $; {\em repo} $\gets \emptyset$\label{alg:line2}
  \While{{\bf true}}
       \If{$\texttt{Empty}(${\em frontier}$)$ and $\texttt{Empty}(${\em repo}$)$\label{alg:line4}}
          \State \Return {\em failure}
       \EndIf
       \If{$\texttt{Empty}(${\em frontier}$)$}
            \State {\em frontier}.put({\em repo}.get())
            \State {\em explored} $\gets \emptyset$  \label{alg:line8}
       \EndIf
       \State {\em node} $\gets $ {\em frontier}.get()
       \If{$node.level = p$ \label{alg:line10}}
            \State \Return $\argmin_{\mathbf{x}\in\mathbb{R}^n} \| \mathbf{Y}_{\ccalI} - \ccalO_{\ccalI}\mathbf{x} \|_2, node.\ccalI $ \label{alg:line11}
       \EndIf
       \State {\em explored}.add($node$)
       \ForAll{$id \in [0,\,1]$ \label{alg:line13}}
         \State $child \gets \texttt{GetChild}(node, id, \mathbf{Y}, \ccalO, \xoverline{w}_i, \epsilon)$
         \If{$child.residual = 0$}
              \If{$child.level-|child.\ccalI| > \s$ \label{alg:line16}}
                  \State {\bf continue}
              \EndIf
             { \If{$child \in $ {\em frontier} {\bf{or}}  {\em explored} \label{alg:line18}}
                  \State {\em repo}.put($child$)
                  \State {\bf continue } \label{alg:line20}
              \Else
                  \State {\em frontier}.put($child$) \label{alg:line22}
              \EndIf}
              %% %% \If{$child \not\in frontier$ \label{alg:line21}}
              %% %%     \State $frontier.put(child)$ \label{alg:line22}
              %% %% \Else
              %% %%      \State $repo.put(child)$ \label{alg:line28}
              %% %% \EndIf
           \EndIf
         \EndFor
         \EndWhile
         \vspace*{0.8em}
\end{algorithmic}
\end{algorithm}

The proposed search algorithm, illustrated in Alg.~\ref{alg:sse}, generates the true attack assignment and estimates the state. It starts by initializing the artificial root and the three node sets, {\em frontier}, {\em explored}, and {\em repo} [lines \ref{alg:line1}-\ref{alg:line2}]. The first set is the priority queue {\em frontier}, which contains nodes that their parents have been expanded and they themselves  are eligible for expansion but have not yet been selected for expansion. The priority is based on the {lexicographical order} between nodes introduced before.  The second set is called {\em explored}, which stores nodes that have been expanded. Storing those nodes avoids repeated search. Finally, the third set is a priority queue {\em repo} which is a repository for nodes for which (i) an equivalent node is in queue {\em frontier}, or (ii) an equivalent node has been expanded and is in the set {\em explored}.
%\footnote{The sets {\em frontier} and {\em explored} exist in a typical search algorithm, where a graph is divided into three disjoint parts: the nodes already explored, the nodes to be explored next (stored in {\em frontier}) and the remaining unexplored nodes. The set {\em frontier} separates the explored nodes from the remaining nodes~\citet{russell2016artificial}. But in this paper, we add another set {\em repo}. {\em frontier} and {\em repo} together separate the explored nodes from unexplored nodes since part of nodes are added to {\em repo} due to the definition of ``equivalence'' of nodes, which will be discussed later.}

After initialization, Alg.~\ref{alg:sse} repeatedly performs the following steps. First, it checks whether the queues {\em frontier} and {\em repo} that contain nodes to be expanded are empty [lines \ref{alg:line4}-\ref{alg:line8}]. If both are empty, then the algorithm has searched the whole graph and it was not able to find a solution. If only {\em frontier} is empty, it is possible that the solution is in the set {\em repo}, {since {\em repo} also contains nodes to be explored next}. In this case, the algorithm picks the node from {\em repo} with the highest priority and puts it in the set {\em frontier}, meanwhile clearing the set {\em explored}. {One can think of selecting a node from {\em repo} as reinitializing the search.} Once a node is included in the set {\em frontier}, Alg.~\ref{alg:sse} proceeds to check and expand it. If this node is at the highest level [lines \ref{alg:line10}-\ref{alg:line11}], the algorithm terminates with a full assignment that satisfies the constraints~\eqref{equ:mina} and~\eqref{equ:minb}; otherwise, this node is expanded [lines \ref{alg:line13}-\ref{alg:line22}]. Note that each node has two possible children nodes. We discuss how to generate a child node later in Alg.~\ref{alg:child}.

Given a generated child node, Alg.~\ref{alg:sse} will discard this node if it violates the constraint~\eqref{equ:mina}. Otherwise, it checks whether the number of attacked sensors in the path leading to this child is larger than $\s$. If this is true, then this node is also discarded {due to the violation of~\eqref{equ:minb}.} If not, then the algorithm searches the set {\em frontier} and {\em explored} for equivalent nodes, i.e., nodes with same values for $value$ and $level$. {If such nodes exist in {\em frontier} or {\em explored}, then the child node is added to {\em repo} and its expansion is delayed [lines \ref{alg:line18}-\ref{alg:line20}]. This node is not discarded because, compared to its equivalent nodes that share the same $level$ and $value$, it is associated with a different partial assignment.} Later we show  that a node with a specific partial assignment exists in {\em frontier} at most once throughout the whole search. {Furthermore, we place the child node into {\em frontier} if no equivalent node exists in these two sets [line~\ref{alg:line22}].}  A simple case study  is shown in Example~\ref{ex:demo}.

\begin{algorithm}[t]
  \caption{\texttt{GetChild} ($parentnode, id, \mathbf{Y}, \ccalO, \xoverline{w}_i, \epsilon$)}
  \label{alg:child}
  \begin{algorithmic}[1]
    \State  $child.level = parentnode.level + 1$ \label{child:line1}
    \State $child.value = id$
    \State  $child.parent = parentnode$
    \If{$id = 0 $ \label{child:line4}}
      \State $child.\ccalI = parentnode.\ccalI \cup \{child.level\}$
      \If{\label{child:min}
        $\min_{\mathbf{x}\in\mathbb{R}^n} \| \mathbf{Y}_{\ccalI} - \ccalO_{\ccalI}\mathbf{x} \|_2 < \xoverline{w}_{\ccalI} + \sqrt{\epsilon}$ \label{child:line5}}
      \State $child.residual =  0 $
        \Else
        \State $child.residual = 1 $ \label{child:line6}
        \EndIf

        \Else \label{child:line7}
    \State $child.\ccalI = parentnode.\ccalI$
    \State $child.residual = parentnode.residual$ \label{child:line9}
    \EndIf
    \end{algorithmic}
\end{algorithm}

Alg.~\ref{alg:child} describes how to generate a child node. The key requirement is that the child node is always one level higher than its parent [lines \ref{child:line1}]. If $child.value = 0$, that is, the corresponding sensor is attack-free, then the algorithm checks whether sensors in the updated set $\ccalI$ of attack-free sensors are consistent meaning that the inequality~\eqref{equ:mina} holds [lines \ref{child:line4}-\ref{child:line6}]. If not, it means the  current partial assignment is wrong in that it incorrectly treats attacked sensors as attack-free. If $child.value = 1$, i.e., the corresponding sensor is under attack, then the set $\ccalI$ remains unchanged and the child node inherits the value $residual$ from its parent [lines \ref{child:line7}-\ref{child:line9}]. This is because the algorithm considers this sensor attacked, so it is not used to reconstruct the state.
  {\begin{exmp}\label{ex:demo}
    \begin{figure}
      \centering
      \includegraphics[width=0.5\linewidth]{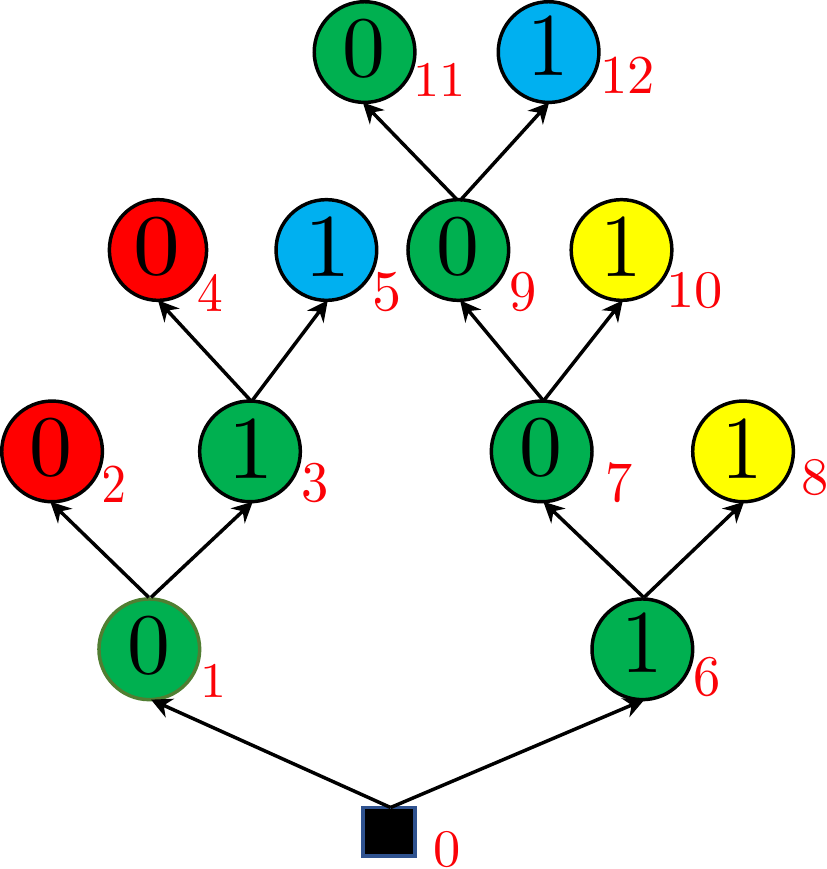}
      \caption{Graphical illustration of Example~\ref{ex:demo}}\label{fig:demo}
    \end{figure}
    {Consider a system with 4 sensors and that only the first sensor is under attack. Assume that this system is 3-sparse observable, which means the observability matrix of each sensor has full rank and we can use only one sensor to reconstruct states. The tree structure in Fig.~\ref{fig:demo} shows how Alg.~\ref{alg:sse} searches the graph in Fig.~\ref{fig:graph}. We differentiate equivalent nodes with different assignment explicitly, thus each path corresponds to one partial/full assignment. Nodes are numbered with indices to their right. Blue, green and yellow nodes are those that are in sets {\em frontier}, {\em explored} and {\em repo} upon termination, respectively; and red nodes are discarded due to the violation of constraints~\eqref{equ:mina} or~\eqref{equ:minb}. The node numbered 11, with tuple $(4, 0, 9, \set{2,3,4}, 0)$, is the first one to reach the highest level, and the path leading to it provides the optimal solution, which is $1,0,0,0$. Thus, only the first sensor is under attack. Table~\ref{tab:demo} shows the node indices in three sets at the beginning of each iteration.}
    \begin{table}[t]
      \centering
      \renewcommand{\arraystretch}{1.2}
       \caption{The evolution of three different sets in Example~\ref{ex:demo} as the Alg.~\ref{alg:sse} proceeds, where ``$-$'' means no nodes exist, and nodes in {\em frontier} and {\em repo} are sorted by priority.}\label{tab:demo}
      \begin{tabular}{clll}
        \hline
        iteration & {\em frontier} & {\em explored} & {\em repo}\\
        \hline
        1 & 0 & $-$ & $-$\\
        2 & 1,6 & 0 & $-$ \\
        3 & 3,6 & 0,1 & $-$\\
        4 & 6,5 & 0,1,3 & $-$\\
        5 & 7,5 & 0,1,3,6 & 8 \\
        6 & 9,5 & 0,1,3,6,7 & 10,8 \\
        7 & 11,12,5 & 0,1,3,6,7,9 & 10,8\\
        8 & 12,5 & 0,1,3,6,7,9,11 & 10,8\\
        \hline
      \end{tabular}
    \end{table}
  \end{exmp}}

\section{Completeness, Optimality, and Complexity}\label{sec:opt}
In the following results, it will help to view two nodes with the same $value$ and $level$ but distinct $\ccalI$ as different. The reason is that, while such nodes are really equivalent, they are treated differently by the algorithm in the sense that a node is added to {\em repo} instead of being discarded  when another node with same $value$ and $level$ is in queue {\em frontier}. In what follows, we show that Alg.~\ref{alg:sse} is complete and optimal and analyze its complexity.

\subsection{Completeness and Optimality}
%% First, we show the completeness and optimality of Alg.~\ref{alg:sse} when dealing with ``perfect'' model, i.e., $ \xoverline{w}= 0$ and $\epsilon=0$.
\begin{thm}\label{thm:comp}
Let the attacked linear dynamical system in \eqref{equ:sys} be $2 \s$-sparse observable and $\epsilon= \epsilon^*$. Assume also that the number of attacked sensors is less than $\bar{s}$ and that each attack signal satisfies $\norm{\mathbf{e}_{i}} > \left( \frac{2}{\sqrt{1-\Delta_s}}\right) \xoverline{w} + \frac{\sqrt{\epsilon}}{\sqrt{1-\Delta_s}}$. %% Considering the noiseless case, where process noise $ \xoverline{w}= 0$ and numerical solver tolerance $\epsilon=0$.
Then Alg.~\ref{alg:sse} is complete that is, if there exists a solution $(\mathbf{x},\mathbf{b})$ that satisfies the constraints~\eqref{equ:mina} and~\eqref{equ:minb} in Problem \ref{pb:min}, then Alg.~\ref{alg:sse} will find it.
\end{thm}

\begin{pf}
  We show  that Alg.~\ref{alg:sse}  terminates in a finite number of iterations at which point the queue {\em frontier} contains a feasible node that satisfies constraints~\eqref{equ:mina} and~\eqref{equ:minb}. A detailed proof of this result is provided in Appendix~\ref{appendixa}.
\end{pf}

%% %% It's possible that for a specific subpath, Alg. \ref{alg:sse} treats an attacked sensor as an attacked-free one and still the $l_2$ residual norm is zero, such as the resulting linear system is underdetermined or full rank. However, those path having this subpath as their prefix part will terminate before they can each the highest level, since the output should satisfy $\mathcal{I} \subseteq \xoverline{\supp(\mathbf{b}^*)}$.

%\subsection{Optimality}
Next, we show that  Alg.~\ref{alg:sse} is optimal, meaning that  Alg.~\ref{alg:sse} can identify true attacks and makes no mistakes in treating attack-free sensors as attacked.
\begin{thm}\label{thm:opt}
  Let the attacked linear dynamical system in \eqref{equ:sys} be $2\s$-sparse observable and $\epsilon= \epsilon^*$. Assume also that the number of attacked sensors is less than $\bar{s}$ and that each attack signal satisfies $\norm{\mathbf{e}_{i}} > \left( \frac{2}{\sqrt{1-\Delta_s}}\right) \xoverline{w} + \frac{\sqrt{\epsilon}}{\sqrt{1-\Delta_s}}$. %% Considering the noiseless case, where process noise $ \xoverline{w}= 0$ and numerical solver tolerance $\epsilon=0$.
Then, the solution of Problem \ref{pb:min} constructed by  Alg. \ref{alg:sse} is optimal meaning that $\mathbf{b}^* = \mathbf{b}_0$, where $\mathbf{b}_0$ represents  the true attack assignment.
\end{thm}

\begin{pf}
 We show that it is not possible that a suboptimal path can be formed before the optimal one. This is because Alg.~\ref{alg:sse} always prioritizes search in directions of possible attack-free assignments and in the case of mistakes, it uses the nodes in the set {\em repo} to take corrective actions. A detailed proof of this result is provided in Appendix~\ref{appendixb}.
\end{pf}

\subsection{Complexity Analysis}\label{sec:comp}
Given a LTI system in~\eqref{equ:sys} and a true attack assignment, in this section we discuss the complexity of Alg.~\ref{alg:sse} in terms of the number of iterations until termination. Note that at each iteration, a node is selected from the queues {\em frontier} or {\em repo}, thus, we can focus on the number of nodes that have been expanded before termination.

%% In the worst case, Alg.~\ref{alg:sse} will search the whole state space, thus the maximum number of iterations is $O(2^p)$. Nevertheless, this number is too conservative, since a large number of invalid nodes will be discarded and the nodes that follow will not be explored any further.
To simplify the analysis of complexity, consider an ``ideal'' model without noise and solution accuracy. Then the constraint~\eqref{equ:mina} becomes
\begin{align}\label{equ:zero}
  \min_{\mathbf{x} \in \rr{R}{n}} \norm{ \mathbf{Y}_{\ccalI} - \ccalO_{\ccalI} \mathbf{x}} = 0.
\end{align}
The following theorem provides an upper bound on the number of iterations taken to find the true assignment.
{\begin{thm}\label{thm:iter}
  Let the attacked linear dynamical system in \eqref{equ:sys} be $2\s$-sparse observable. Assume also that the true number of sensors that are  under attack is $s$. Let $S = p - 2\s$, where $p$ is the number of sensors and $\s$ is the maximum allowable number of sensors under attack. Then, without considering noise and solution accuracy, {Alg.~\ref{alg:sse} takes at most {$N_{\textup{upper}} = \sum_{i=1}^{S} {s \choose i} {\s+S-s \choose S-i} (\s + S)+p$} iterations to find the true assignment, where ${s \choose i} =0$ if $s < i$.}
\end{thm}}
\begin{pf}
  Alg.~\ref{alg:sse} achieves its  worst performance if the first $s$ sensors are the ones that are under attack. {This is because, in this case, the attacked sensors can be treated as attack-free when the system of equations $\mathbf{Y}_{\ccalI} =  \ccalO_{\ccalI} \mathbf{x} $ (where $\ccalO_{\ccalI}\in \mathbb{R}^{T\cdot |\ccalI|\times n}$) is underdetermined or square, i.e., when $T \cdot |\ccalI| \leq n$, and Alg.~\ref{alg:sse} is biased towards expanding attack-free sensors first.} Recalling  the graph in Fig.~\ref{fig:graph}, this worst scenario corresponds to the case where the first $s$ levels are associated with sensors that are under attack.  In what follows, we focus on this worst-case attack scenario where the first $s$ sensors are under attack.

  Since the system is $2\s$-sparse observable, any observability matrix $\ccalO_{\ccalI}$ corresponding to  $\ccalI$ with $|\ccalI| \geq p-2\s$  has full rank $n$. Assume first that $s \geq S = p - 2\s$. {When $|\ccalI|<S$, the rank of the observability matrix $\ccalO_\ccalI$ can be smaller than $n$.  We further assume that the system of linear equations $\mathbf{Y}_{\ccalI} =  \ccalO_{\ccalI} \mathbf{x} $ is feasible, which could occur when it is underdetermined or square.} Then,
  the algorithm can assign 0's to $S$ sensors out of the first $s$ sensors.
  %% because $\mathbf{x} = \ccalO_{\ccalI}^{-1}\mathbf{Y}_{\ccalI}$ satisfies~\eqref{equ:zero} when $|\ccalI|=S$.
  But after the $S$-th sensor that is assigned 0, if the next sensor is also assigned 0, then the system of linear equations corresponding to those $S+1$ sensors, which are treated as attack-free, becomes inconsistent since it incorrectly treats  attacked sensors as attack-free. Therefore, the algorithm can only assign 1 to all sensors past the $S$ sensors that have been assigned 0 until there are $\s$ 1's in this path. Hence, in this case, the whole searched path has {$\s + S$} nodes. Similarly, if $S-1$ nodes are assigned 0 and $s - (S-1)$ nodes are assigned 1 for the first $s$ sensors, Alg.~\ref{alg:sse} can assign 0 to one more sensor and 1 to another $\s-[s- (S-1)]$ sensors. Following this logic, we get that if $i$ nodes are assigned 0 among the first $s$ sensors, where {$1\leq i \leq S$}, Alg.~\ref{alg:sse} can assign 0 to another $S-i$ sensors and 1 to $\s - (s - i) $ sensors. We conclude that our algorithm can only explore $i + (S-i) + (s-i) + [\s - (s-i)] = \s+S $ nodes on one path if this path is not feasible. Each assignment to the first $S+\s$ sensors corresponds to one path. In the worst-case scenario considered here, the path leading to the optimal solution has the lowest priority in the sense that a node at a specific level $l \leq S$ on an infeasible  path has higher priority  than the node with same $value$ and $level$ on the optimal path. Hence, all infeasible paths will be searched before the algorithm terminates, and we denote by $\Pi$ the set of these paths.\footnote{Note that infeasible paths can only be partially searched. As shown in Theorem~\ref{thm:opt}, it is not possible that infeasible paths can be fully searched before the optimal path.} Therefore, the search algorithm searches at most {$\sum_{i=1}^{S} {s \choose i} {\s+S-s \choose S-i} (\s + S)$} nodes corresponding to infeasible paths. When no nodes are assigned 0 at the first $s$ levels, Alg.~\ref{alg:sse} will assign 0's to the remaining {$p-s$} sensors, which is exactly the number of nodes needed to be required corresponding to the worst-case attack scenario. Because Alg.~\ref{alg:sse} explores one node per iteration, we get the upper bound on the number of iterations taken to find the true assignment.

  When $s< S$, the situation is less complex since the path corresponding to the optimal solution is inside $\Pi$, that is, the path where the first $s$ nodes are assigned 1 is inside $\Pi$, so there is no need to explore the graph until the first $S$ nodes are assigned 1, completing the proof.
\end{pf}
\begin{rem}
  {%% The upper bound $\sum_{i=1}^{S} {s \choose i} {\s+S-s \choose S-i} (\s + S)+p$ can be rewritten as $\sum_{i=1}^{p-2\s} {s \choose i} {p-\s-s \choose p-2\s-i} (p-\s)+p$ due to $S=p-2\s$. When $s$ and $p$ are fixed, this upper bound decreases as $\s$ increases. Consider $\s' =  \s + 1$, and the corresponding $N_{\text{upper}}$ and $N'_{\text{upper}}$, for any term ${s \choose i} {p-\s'-s \choose p-2\s'-i} (p-\s')$ in $N'_{\text{upper}}$, there always exists a term ${s \choose i} {p-\s-s \choose p-2\s-i} (p-\s)$ that is larger than it in $N_{\text{upper}}$.
    {(Ideal worst case): {The worst-case complexity in Theorem~\ref{thm:iter} assumes  that if  $\mathbf{Y}_{\ccalI} =  \ccalO_{\ccalI} \mathbf{x}$ is underdetermined or square, i.e., if  $T \cdot |\ccalI|\leq n$, then feasible solutions exist. However, in practice, it is possible that no solution exists.  Furthermore, as the index set $\ccalI$ grows, the system of linear equations becomes overdetermined ($T \cdot |\ccalI| > n$), which is almost always inconsistent and thus has no solutions. In this case,  the wrong assignment will terminate much sooner, as demonstrated in Section~\ref{sec:noiseless}. Therefore, we refer to the worst case in the proof  as the ``ideal'' worst case. Whether the ideal worst case might occur is application-dependent, and relies on the dynamical system,  the number of measurements $T$ and injected attack signals.}}}
\end{rem}

\section{Experimental Results}\label{sec:sim}
In this section, we present several test cases for values of $p$ and $n$, implemented using Python 3.6.3 on a computer with 2.3 GHz Intel Core i5 and 8G RAM, that illustrate the correctness and  efficiency of the proposed algorithm {for large-scale estimation problems}. To validate our method, we compare  with the Mixed Integer Quadratically Constrained Programming (MIQCP) method in~\citet{winston2004operations} and with the solver IMHOTEP-SMT.\footnote{IMHOTEP-SMT source code  in Matlab  can be found at \url{http://nesl.github.io/Imhotep-smt/index.html}.} The formulation of the MIQCP problem  takes the form
  \begin{align}\label{equ:miqcp}
    \min_{(\mathbf{x},\mathbf{b})\in \rr{R}{n} \times \rr{B}{p}}  & \;  \sum_{i=1}^{p} b_i \\
    s.t. \quad  \; \norm{\mathbf{Y}_i - & \ccalO_i \mathbf{x}}  \leq M b_i + \xoverline{w}_i + \sqrt{\epsilon}, \;\forall i \in \set{1,\ldots, p}. \nonumber
  \end{align}
  Since $M\in {\rr{R}{}}$ is a very big number, the constraints in \eqref{equ:miqcp} are called  Big-M constraints, which are used to model the binary activation/deactivation of constraints defined over real decision variables. The performance of MIQCP is sensitive to the value of $M$. Given different $M$, the resulting attack assignment  may vary significantly. %% We have to carefully tune $M$ to obtain a satisfied solution, but when we do not know the actual attack, there is little information to leverage when tuning $M$.
  Note that we can get $\mathbf{x}$ and $\mathbf{b}$ simultaneously by solving the MIQCP, but in our experiments, using the commercial solver Gurobi \cite{gurobi}, we observed that selecting a very large $M$ can provide a feasible assignment $\mathbf{b}$ but a bad state $\mathbf{x}$. This is because using large $M$ places more emphasis on the optimization over the binary variables $\mathbf{b}$ rather than the real variables $\mathbf{x}$. To overcome this issue, we implement MIQCP in two steps. First, we solve the MIQCP problem~\eqref{equ:miqcp} for large $M$ to obtain a feasible assignment $\mathbf{b}$ and then use this assignment to solve the unconstrained least square problem
\begin{align*}
  \min_{\mathbf{x} \in \rr{R}{n}} \norm{\mathbf{Y}_{\ccalI} - \ccalO_{\ccalI} \mathbf{x}},
\end{align*}
for the state $\mathbf{x}$, where the set $\ccalI$ is the set of attack-free sensors predicted by the solution of  the MIQCP. On the other hand, IMHOTEP-SMT is mainly designed to find a feasible solution, and it can be used to find the optimal solution in the noiseless case if executed repeatedly for different values of $\bar{s}$ by performing a binary search over $\bar{s}$, checking feasibility of the system of equations in~\eqref{equ:mina}, decreasing $\bar{s}$ if these equations are feasible, and repeating this process until the constraint~\eqref{equ:mina} is violated or until $\bar{s}=0$. However, IMHOTEP-SMT needs to search all possible combinations of truth assignments to all sensors to ensure that a problem is infeasible, thus selecting a feasible but less conservative $\bar{s}$ for IMHOTEP-SMT is not easy. Hence, we run IMHOTEP-SMT once until the first feasible solution is found.

We randomly construct sparse matrices $A$ and $C$ for various parameters $p$ and $n$ with entries in the interval $[0,1]$. Furthermore, the initial state and the set of the attacked sensors and corresponding attack signals are also randomly generated, as in~\citet{fawzi2014secure,shoukry2017secure,pajic2015attack,pajic2017attack}. {Specifically, given an LTI system, we  adopt two schemes to select the set of attacked sensors~\cite{park2017security}. The first scheme attacks the first $s$ sensors in the graph in Fig.~\ref{fig:graph}, which we refer to as the greedy attack, and the second scheme randomly selects $s$ sensors to attack. For the random attack scheme, we randomly generate $m$ different {true attack assignments}, and attack signals that satisfy inequality~\eqref{equ:attack} over $T=n$ time instants for each attack assignment. {Specifically, the attack signals are generated by first sampling a vector which follows a standard normal distribution, then normalizing it, and finally multiplying by a number which captures the magnitude, as in~\citet{shoukry2017secure}.} For the greedy scheme, we randomly generate $m$ such attack signals. The solution  accuracy is set to $\epsilon=10^{-5}$. The quantity $M$ in MIQCP is set to be $10^8$.} We monitor the execution time and relative estimation error of each method, defined as $\frac{\norm{\mathbf{x}_0 - \mathbf{x}^*}}{\norm{\mathbf{x}_0}}$. Note that $\mathbf{x}_0$ is the true state, while with a slight abuse of notation, $\mathbf{x}^*$ is the output of each method.

{\subsection{Optimality and Complexity}
In this section, we validate the optimality and complexity of the proposed algorithm. For this, we need to get the maximum allowable number of attacked sensors $\bar{s}$ by increasing $\bar{s}$ from 1 to $\ceil*{p/2}-1$ and checking the observability matrix after excluding any $2\bar{s}$ sensors. Since this problem is combinatorial, we consider small systems with $n=p=10$ for which $\bar{s}$ and the resulting $\Delta_s$ in inequality~\eqref{equ:attack} are computationally tractable. Then, we validate optimality by generating attack signals that satisfy~\eqref{equ:attack}, and complexity by comparing the number of iterations taken with the upper bound.}

\begin{table}[t]
  \caption{Properties of LTI systems}\label{tab:lti}
  \centering
\renewcommand{\arraystretch}{1.1}
  \begin{tabular}{ccccccc}
    \toprule
    \multirow{2}{*}{$\s$} & \multicolumn{4}{c}{$\Delta_s$} &  \multirow{2}{*}{$\frac{2}{\sqrt{1-\Delta_s}}$}  &  \multirow{2}{*}{$N_{\text{upper}}$}\\
    \cline{2-5}
    & mean & std & min & max & & \\
    2 & 0.985 & 0.024 & 0.894 & 0.999994 & 16.341 & 226\\
    3 & 0.986 & 0.021 & 0.910 & 0.999890 & 17.202 & 248\\
    4 & 0.999 & 0.002 & 0.994 & 0.999994 & 55.541 & 94\\
    \midrule
  \end{tabular}
\end{table}

\begin{figure}[t]
  \centering
  \includegraphics[width=0.85\linewidth]{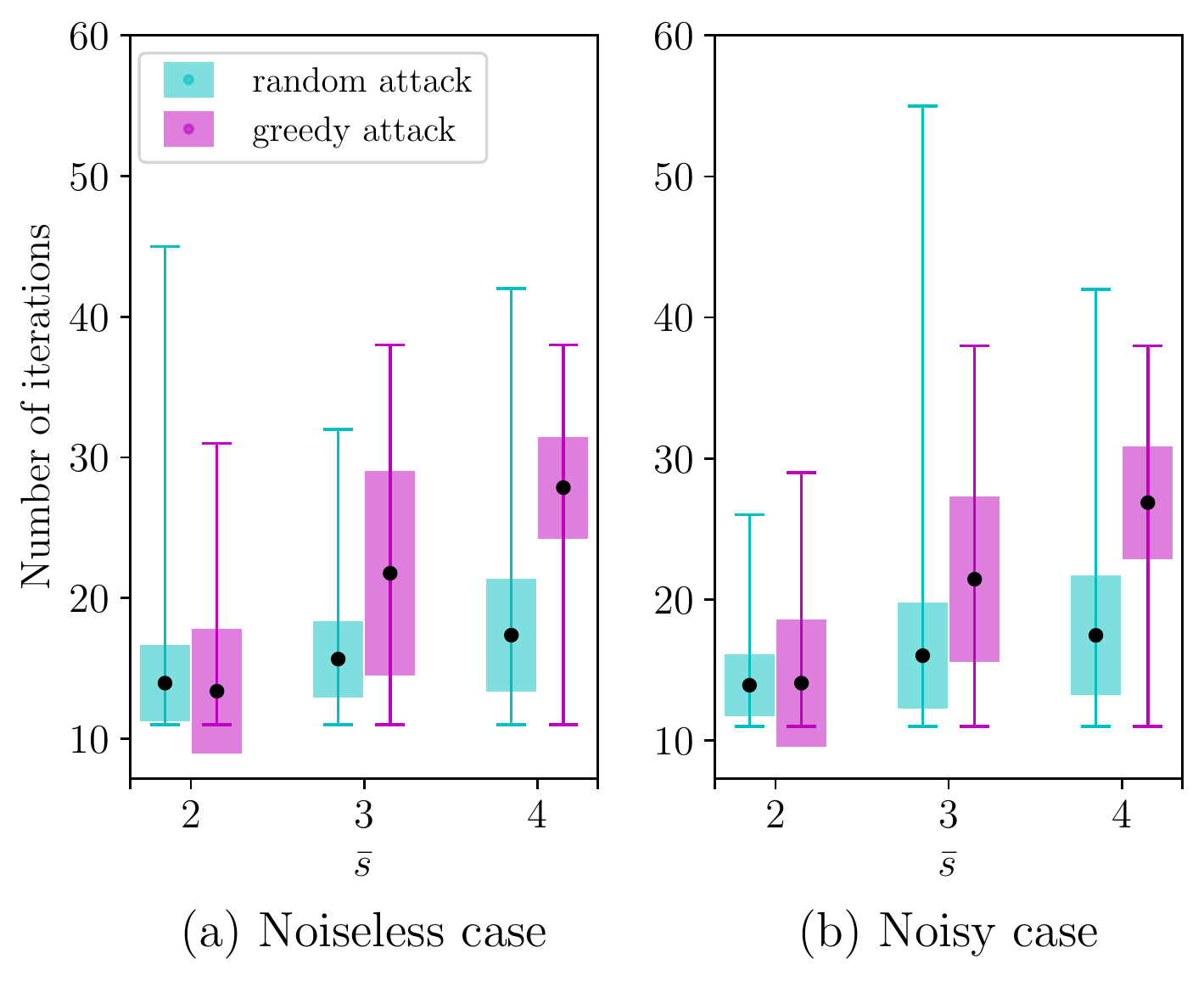}
  \caption{Secure state estimations for small systems}\label{fig:small}
\end{figure}

{Specifically, we consider both the noiseless and noisy cases. For each case, we randomly generate three sets of $2\bar{s}$-sparse-observable and small-scale systems, 50 systems per set,  with $p=n=10$, such that each set corresponds to systems with $\bar{s}=2, 3, 4$, respectively. Then, for each system,  we set the number of attacked sensors $s$ to be $\bar{s}$. Table~\ref{tab:lti} records the statistics, including the mean, standard deviation, and minimum and maximum value of $\Delta_s$ with respect to the maximum allowable number of attacked sensors $\s$, computed over the 50 systems per set. {The column $\frac{2}{\sqrt{1-\Delta_s}}$ in Table~\ref{tab:lti} shows the mean value and represents the average power of the attack signal compared to that of the noise.} The last column shows the theoretical upper bound on the maximum number of iterations that Alg.~\ref{alg:sse} will take given $p$, $\s$ and $s=\s$. For the noiseless and noisy cases, we compare the performance of Alg.~\ref{alg:sse} for different values of $\s$ and for the two different attack schemes. From our simulations, we observe that the true attack assignment can be identified every time for both the noiseless and noisy cases. Furthermore, in~Fig.~\ref{fig:small} we report the statistics on the number of iterations taken over 1000 trials per case. It can be seen that, on average, identifying the true attack assignment for the greedy attack requires more iterations compared to the random attack. {For   small-scale systems, we found that when the index set $\ccalI$ contains only one sensor that is also under attack, it is often the case that the square matrix $\ccalO_\ccalI$, with $T=n$, is non-singular and  thus the system of linear equations $\mathbf{Y}_{\ccalI} =  \ccalO_{\ccalI} \mathbf{x}$ is consistent. In this case, the sensor will be mistakenly identified as attack-free.  Thus, more iterations are required,  which verifies the complexity analysis for the ideal worst case presented in Section~\ref{sec:comp}. }
  Notably, the number of iterations of Alg.~\ref{alg:sse} is lower than the  upper bound regardless of the presence or absence of the noise and the attack scheme. {This is because when the index set $\ccalI$ includes more than one sensors and at least one of those sensors  is under attack, the system of linear equations is overdetermined and inconsistent, so that the wrong assignment will terminate immediately.}}
\begin{figure}[t]
  \centering
  \includegraphics[width=0.85\linewidth]{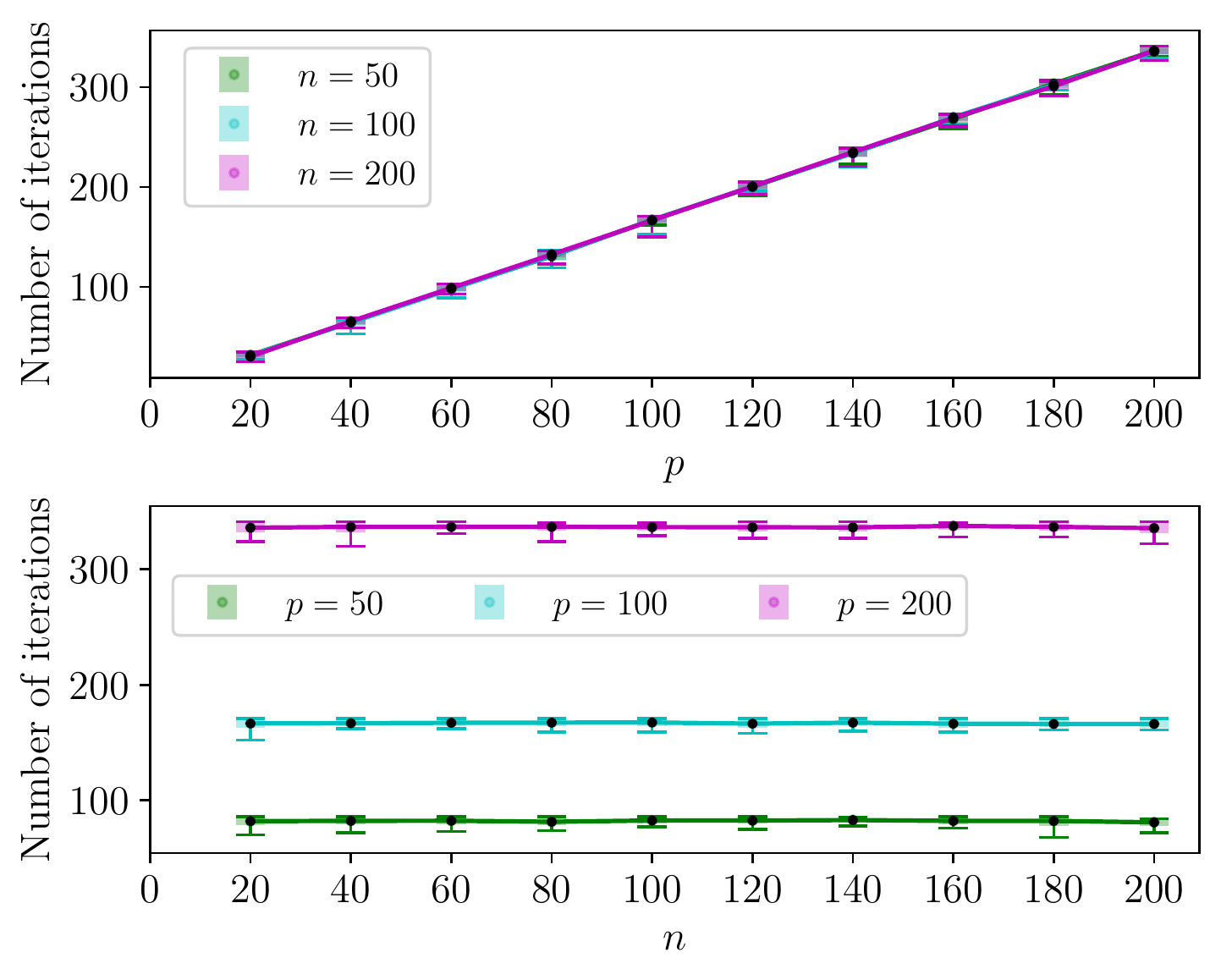}
  \caption{Number of iterations versus scales of systems}
  \label{fig:complexity}
\end{figure}
%% \begin{figure*}[t]
%%   \centering
%%   \subfigure[$s/p=10\%$]{
%%     \label{fig:noiseless_1}
%%     \includegraphics[width=0.3\linewidth]{figure/noiseless_1.pdf}}
%%    \subfigure[$s/p=20\%$]{
%%     \label{fig:nioseless_2}
%%     \includegraphics[width=0.3\linewidth]{figure/noiseless_2.pdf}}
%%    \subfigure[$s/p=30\%$]{
%%     \label{fig:noiseless_3}
%%     \includegraphics[width=0.3\linewidth]{figure/noiseless_3.pdf}}
%%   \caption{Runtimes versus scales of systems for the noiseless case}\label{fig:noiseless}
%% \end{figure*}

\begin{figure}[t]
  \centering
  \includegraphics[width=0.85\linewidth]{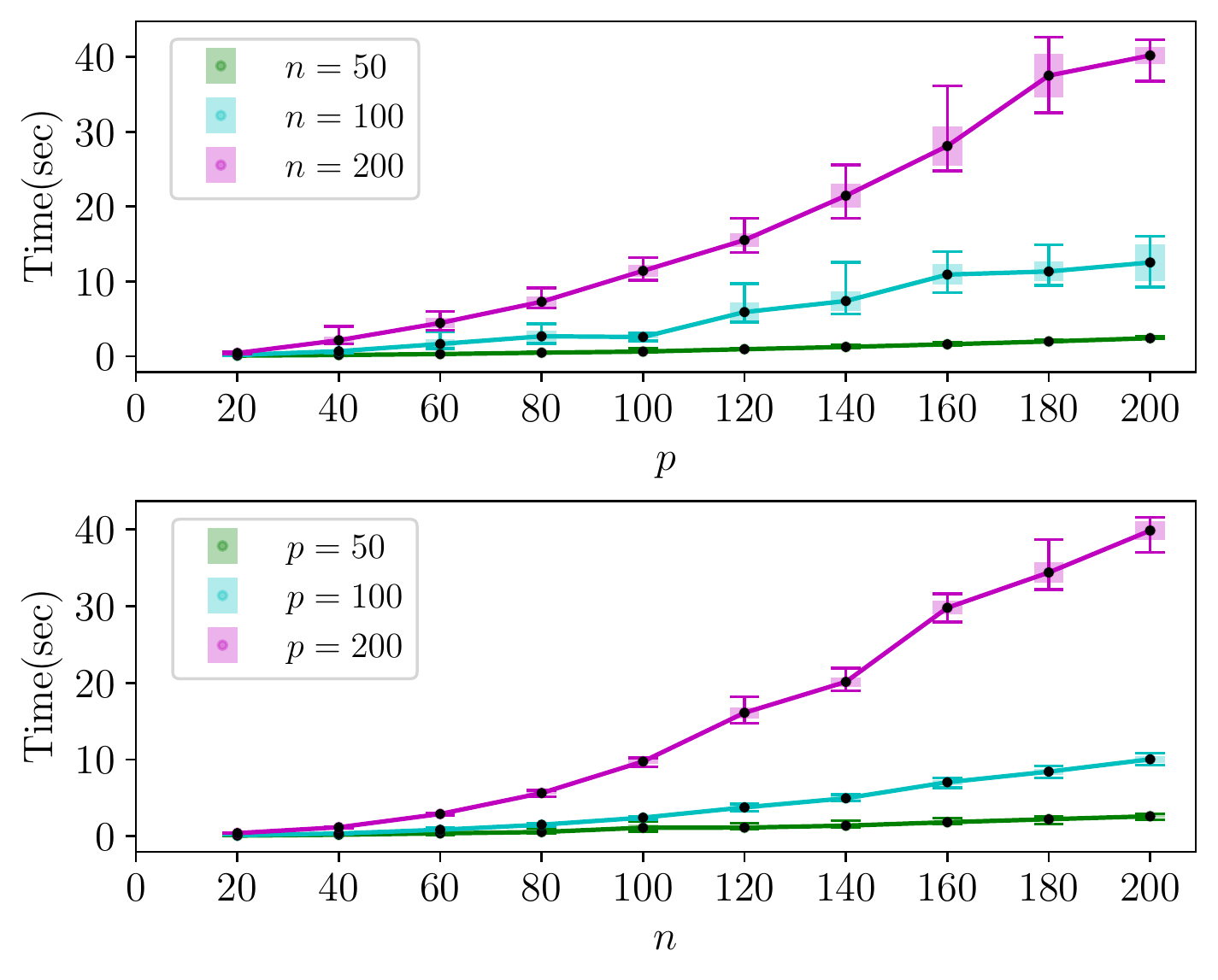}
  \caption{Runtimes versus scales of systems}
  \label{fig:complexity_time}
\end{figure}

{{Next,  we demonstrate that the upper bound on the number of iterations given in Thm.~\ref{thm:iter} is independent of the number of states. For this, we assume no process or measurement noise. Specifically, we first fix $n$ at values 50, 100, 200, respectively, and then gradually increase the number of sensors up to 200. Then, we fix $p$ at values 50, 100, 200, respectively, and gradually increase the number of states up to 200. For each combination of $p$ and $n$, we randomly generate $5$ LTI systems.  For each LTI system, {the number of attacked sensors is chosen to equal  $30\%$ of the total number of sensors,} and the corresponding attack assignment is generated according to the random attack  scheme with $m=5$.  The number of iterations of Alg.~\ref{alg:sse} averaged over 25 trials are shown in Fig.~\ref{fig:complexity}. It can be seen that the number of iterations increases as $p$ grows but it is not affected by $n$.} We show the execution time in Fig.~\ref{fig:complexity_time}. Although the number of iterations of Alg.~\ref{alg:sse} for the same $p$ does not change significantly with $n$, the runtime of the algorithm for the same $p$ increases as $n$ increases. This is due to the time required to solve the minimization problem in line~\ref{child:line5}, Agl.~\ref{alg:child}.
}

\begin{figure}[t]
  \centering
    \includegraphics[width=0.85\linewidth]{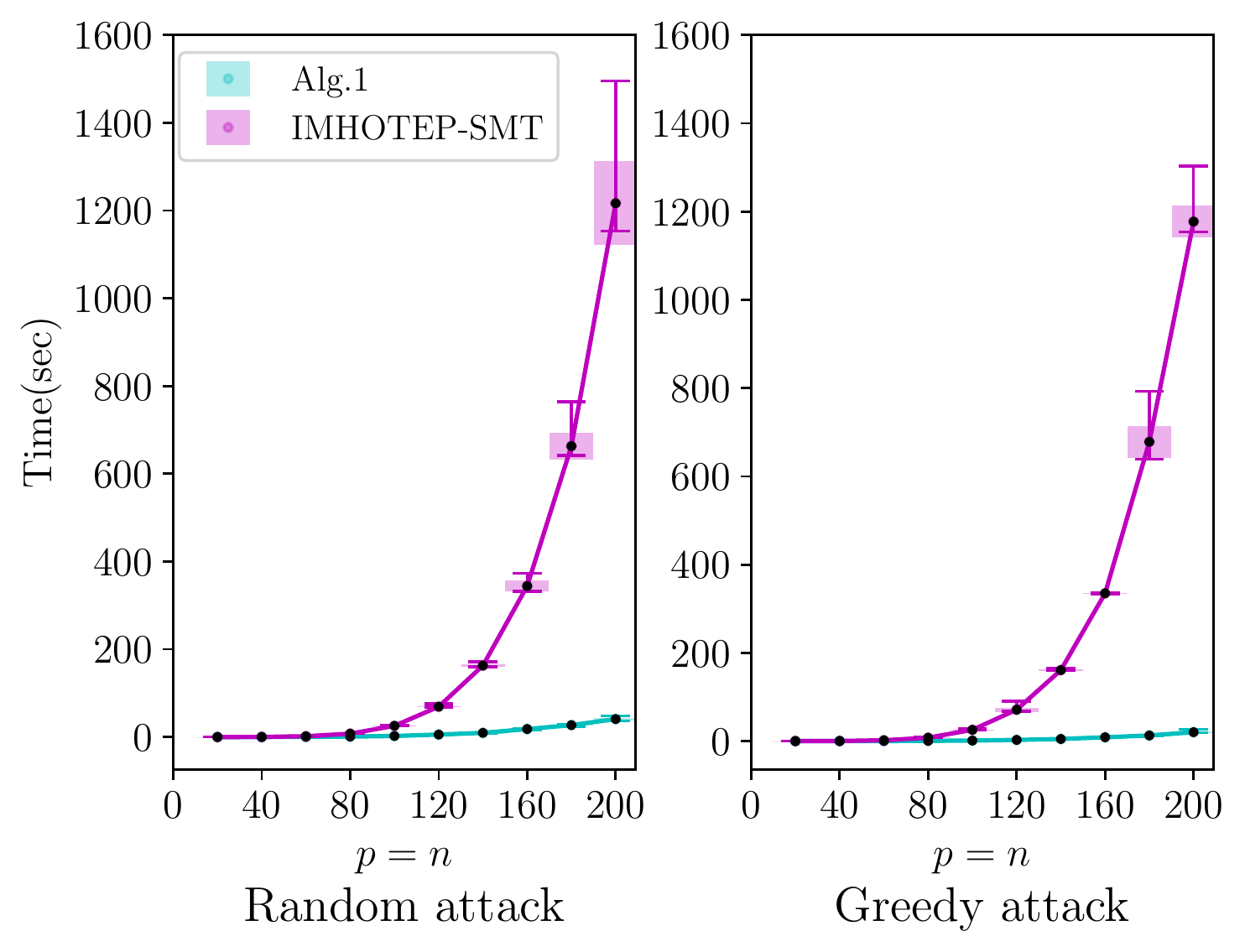}
  \caption{Runtimes versus scales of systems for the noiseless case}\label{fig:noiseless}
\end{figure}
\subsection{Secure State Estimation without Noise}\label{sec:noiseless}
%% \begin{figure*}[t]
%%   \centering
%%   \subfigure[$s/p=10\%$]{
%%     \label{fig:noisy_1}
%%     \includegraphics[width=0.3\linewidth]{figure/noisy_1.pdf}}
%%    \subfigure[$s/p=20\%$]{
%%     \label{fig:niosy_2}
%%     \includegraphics[width=0.3\linewidth]{figure/noisy_2.pdf}}
%%    \subfigure[$s/p=30\%$]{
%%     \label{fig:noisy_3}
%%     \includegraphics[width=0.3\linewidth]{figure/noisy_3.pdf}}
%%   \caption{Runtimes versus scales of systems for the noisy case}\label{fig:noisy}
%% \end{figure*}

\begin{figure}[t]
  \centering
  \includegraphics[width=0.85\linewidth]{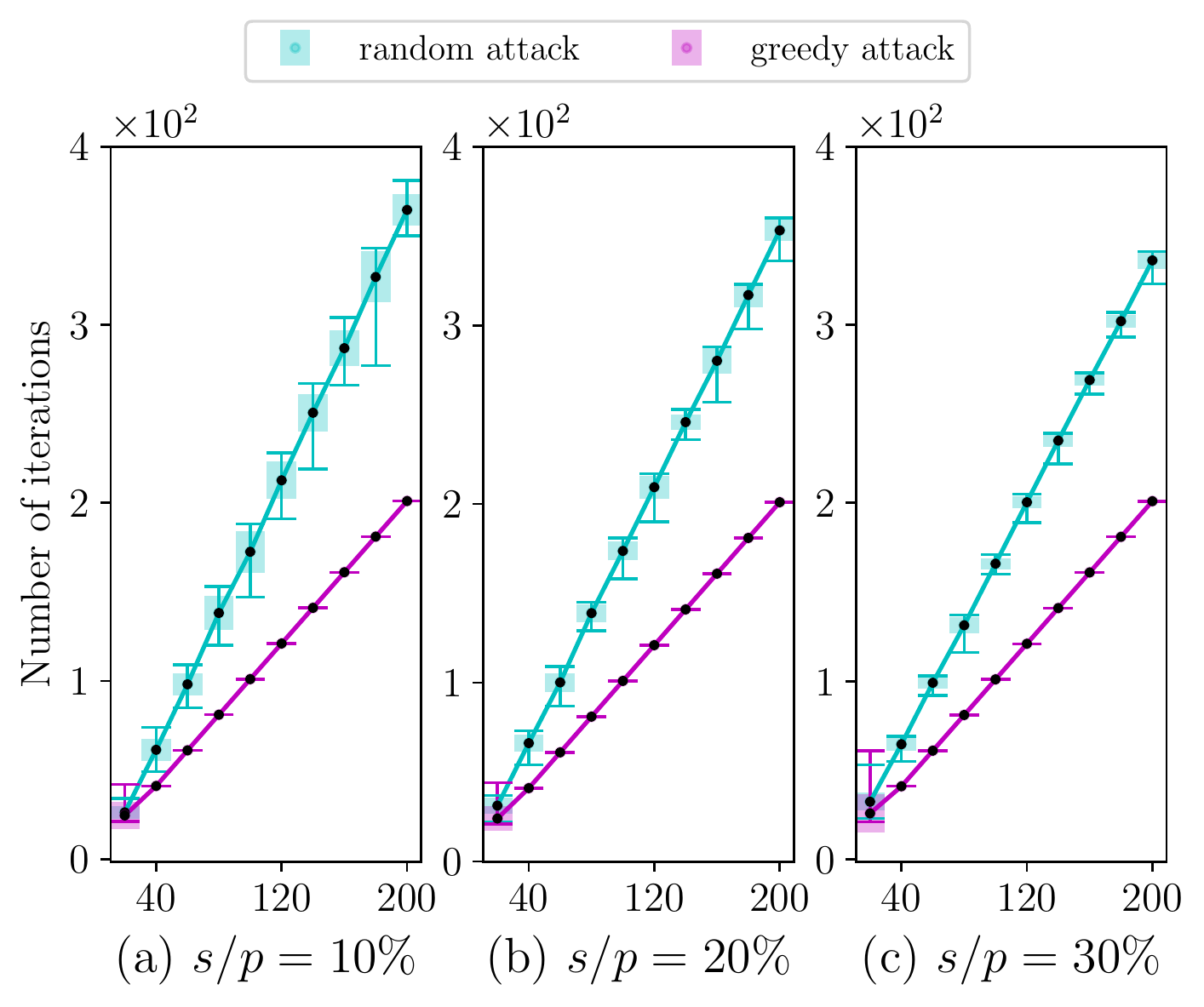}
    \caption{Number of iterations taken versus system scales}\label{fig:iter_noiseless}
\end{figure}

In this section, we consider the noiseless scenario where the number of states $n$ and the number of states $p$ are the same, both varying from 20 to 200. {For each combination of $p$ and $n$, we randomly generate $5$ LTI systems. For each LTI system, {the number of attacked sensors is chosen to equal $10\%, 20\% , 30\%$ of the total number of sensors, respectively,} and the corresponding attack support is generated according to the greedy and random attack schemes, respectively, with $m=5$.} Since it is computationally expensive to obtain $\s$ for large-scale systems, line~\ref{alg:line16} in~Alg.~\ref{alg:sse} is replaced by  $child.level-|child.\ccalI| \geq \ceil*{p/2}$.  The solutions returned by Alg.~\ref{alg:sse} and IMHOTEP-SMT are identical to the true attack assignment in each trial, and their estimation error is close to 0. However, this is not the case for MIQCP. The statistics of the runtimes averaged over 25 trials for $s/p=30\%$ are shown in Fig.~\ref{fig:noiseless}. {The results  for  $s/p=10\%, 20\%$ are similar to those shown in Fig.~\ref{fig:noiseless} and thus, we omit them due to space limitations.}
When the number of states and sensors are small, execution times are comparable for Alg.~\ref{alg:sse} and IMHOTEP-SMT. However, as the system scale becomes large, Alg.~\ref{alg:sse} significantly outperforms IMHOTEP-SMT in terms of runtime. {Specifically, for the random attack  scheme and $n=p=200$, Alg.~\ref{alg:sse} requires 57.7$s$, 49.7$s$ and $41.0s$ on average to find the truth assignment when $s/p=10\%, 20\%, 30\%$, respectively, while IMHOTEP-SMT requires 691.1$s$, 992.2$s$ and 1216.7$s$ and these times grow  much faster than Alg.~\ref{alg:sse}. Fig.~\ref{fig:iter_noiseless} shows  the number of iterations required by Alg.~\ref{alg:sse}, which grow almost linearly with  $p$ and $n$. Note that fewer iterations are required in the greedy attack case, contrary to the result in Thm.~\ref{thm:iter} where the ideal worst case is analyzed. {For large-scale systems, we found that when the index set $\ccalI$ contains only one sensor and this sensor is also under attack, it is more likely  that the square matrix $\ccalO_\ccalI$ is singular so that the  system of linear equations $\mathbf{Y}_{\ccalI} =  \ccalO_{\ccalI} \mathbf{x}$ is  inconsistent.  Therefore, the ideal worst case does not take place and the inconsistency causes the wrong partial assignment to terminate early.} Although it is computationally almost impossible to calculate the maximum allowable number of attacked sensors $\s$ and the upper bound $N_{\text{upper}}$ for large systems, we here provide the upper bound $N_{\text{upper}}$ when $\s=1, 2, 3, 97, 98, 99$, respectively, for $p=200, s/p=10\%$. Approximately, these bounds  are $39801, 3.86\times10^6, 2.47\times10^8, 1.08\times10^{11}, 2.55\times 10^8, 1.83\times10^5$, respectively. Calculating $N_{\text{upper}}$ for other $\s$ is not achievable within reasonable amount of time. Observe that in Fig.~\ref{fig:iter_noiseless}(a) the maximum number of iterations required by Alg.~\ref{alg:sse} is lower than 400, which makes our algorithm at least two orders of magnitude faster than the above upper bound.}

%% To see the actual number of iterations taken compared with the upper bound,
\begin{table}[t]
  \caption{Misidentification ratio of MIQCP for the noiseless case}\label{tab:miqcp_noiseless}
  \centering
  \renewcommand{\arraystretch}{0.9}
  \setlength{\tabcolsep}{2.8pt}
    \begin{tabular}{ccccccccccc}
      \toprule
      $s/p$ & 20 & 40 &60& 80 & 100 & 120 & 140 & 160 & 180 & 200\\
      \midrule
    0.1 &  0.24 & 0.0 &  0.2 & 0.64 & 0.24 & 0.52 & 0.88 & 0.52 & 0.6 & 0.8 \\
    0.2 &  0.32& 0.04& 0.24& 0.56& 0.2& 0.6& 0.8& 0.48& 0.48& 0.8\\
    0.3 & 0.44& 0.08& 0.24& 0.56& 0.28& 0.52& 0.8& 0.52& 0.56& 0.72\\
    \midrule
    0.1 & 0.2& 0.0& 0.2& 0.56& 0.24& 0.56& 0.8& 0.44& 0.56& 0.8\\
    0.2 & 0.24& 0.08& 0.24& 0.6& 0.24& 0.52& 0.8& 0.48& 0.52& 0.92\\
    0.3 & 0.6& 0.12& 0.32& 0.44& 0.2& 0.56& 0.8& 0.48& 0.64& 0.84 \\
    \bottomrule
  \end{tabular}
\end{table}

{As for MIQCP, Table~\ref{tab:miqcp_noiseless} shows the misidentification ratio between  the number of trials where the truth assignment cannot be identified divided by the total number of  25 trials, where the first 3 rows show the results for the random attack scheme and the last 3 rows for the greedy scheme.} It can be seen that the performance of MIQCP is heavily affected by the selection of big $M$ and, in simulation, it behaves poorly in terms of estimation accuracy. {Note that, in practice, tuning $M$ does not help recover the true attack assignment, as suggested by extensive numerical simulations we have conducted. Moreover, there is no way to justify the selection of $M$ if the truth assignment is not known beforehand.}

Since IMHOTEP-SMT has been already compared to the event-triggered projected gradient descent method (ETPG) \cite{shoukry2016event} in terms of execution time and relative estimation error and has been shown to exhibit better performance, we did not compare with this method here. Moreover, we did not compare with methods that relax $\ell_0$-based formulations to convex problems since they lack correctness guarantees.

\subsection{Secure State Estimation with Noise}\label{sec:noisy}

In this part, we consider the noisy case by following the same simulation procedure as in the noiseless case. Note that it is computationally expensive to calculate the term $\Delta_s$  when $p$ and $n$ are large. Instead, we estimate $\Delta_s$ by simulating the process and measurement noises subject  to truncated normal distributions. Then, we select attack signals with relatively large strengths. Since the selected attack vectors can be much stronger than the noise, which is a strategy that a true attacker would not typically follow, we numerically tune the strength of the attack vector by running a large number of trials and checking if IMHOTEP-SMT is able to correctly identify the truly attacked sensors in all trials. If yes, we reduce the strength of the attack vector and repeat this process until our algorithm fails to identify the true attack assignment in all trials.  When this happens, we select values for the attack vector from the last successful trial. Then, Alg.~\ref{alg:sse} is applied to the same attacked dynamical system and corresponding measurements. It turns out that Alg.~\ref{alg:sse} and IMHOTEP-SMT can identify the truth assignment, which validates the workaround above when $n$ and $p$ are large.

\begin{figure}[t]
  \centering
     \label{fig:noisy_3}
   \includegraphics[width=0.85\linewidth]{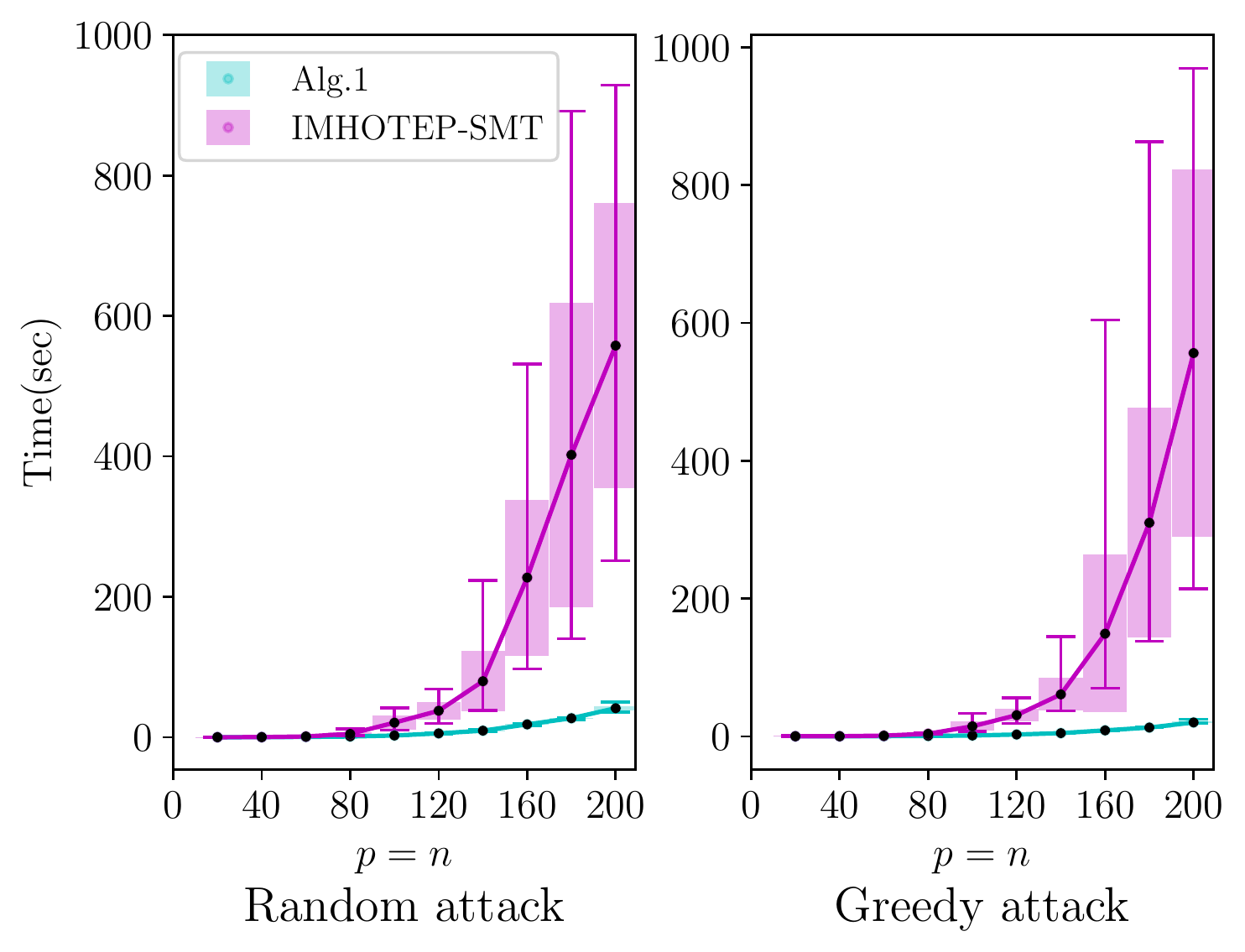}
  \caption{Runtimes versus scales of systems for the noisy case}\label{fig:noisy}
\end{figure}

\begin{figure}
  \centering
  \includegraphics[width=0.85\linewidth]{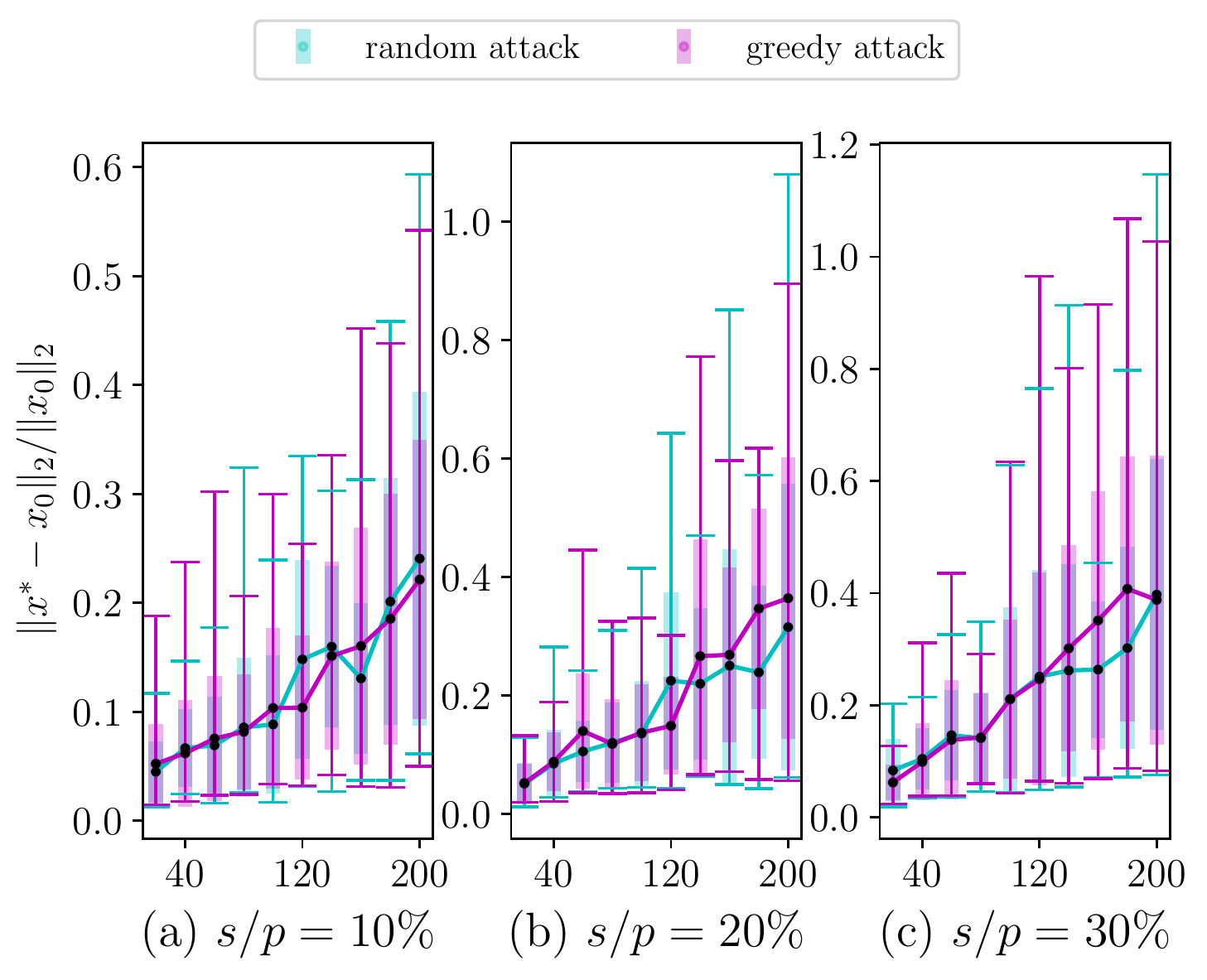}
  \caption{Relative estimation errors versus scales of systems}\label{fig:error}
\end{figure}

{{As in Section~\ref{sec:noiseless}, we monitor the execution time averaged over 25 trials. The results for $s/p=30\%$ are shown in Fig.~\ref{fig:noisy}. The results   for $s/p=10\%, 20\%$ are similar to those shown in Fig.~\ref{fig:noisy} and thus, we omit them due to space limitations.} We observe that, Alg.~\ref{alg:sse} outperforms IMHOTEP-SMT in terms of runtime. Specifically, for $p=n=200$, Alg.~\ref{alg:sse} requires  55.1$s$, 47.5$s$ and 41.3$s$ on average to return the solution for $s/p=10\%, 20\%$ and $30\%$, respectively, while IMHOTEP-SMT requires $275.4s$, $526.1s$ and $557.7s$ and this runtime grows much faster than Alg.~\ref{alg:sse}, as in the noiseless case.  The relative estimation errors of Alg.~\ref{alg:sse} are shown in Fig.~\ref{fig:error}. Furthermore, the misidentification ratio of MIQCP is similar to Table~\ref{tab:miqcp_noiseless} and we omit showing it due to space limitations.}


\subsection{Discussion}
Our numerical results for the noiseless and noisy case studies in Sections~\ref{sec:noiseless} and~\ref{sec:noisy}, show that our method outperforms  IMHOTEP-SMT in terms of runtime. The reason is that, in IMHOTEP-SMT, the optimization solver only checks the correctness of a full truth assignment for all sensors provided by the SMT solver. However, our method also checks partial assignments, which allows a search to terminate early if an assignment is infeasible. Second, in IMHOTEP-SMT, certificates that exploit the geometry of the problem are explicitly generated to serve as heuristics targeting a smaller set of sensors where at least one sensor is attacked, when a full set of  truth assignments is not successful. The idea is that the affine half-spaces $\mathbf{Y}_\ccalI = \ccalO_{\ccalI} \mathbf{x}$ corresponding to the attack-free sensors should intersect, and a small set of affine subspaces failing to intersect means there exists at least one attacked sensor that is incorrectly identified  as attack-free. {In our algorithm, prioritizing search along paths with more attack-free sensors reduces the dimension of the kernel space of the equations $\mathbf{Y}_\ccalI = \ccalO_{\ccalI} \mathbf{x}$, which makes the search less tolerant to mistakes. Therefore, assigning higher priority to nodes in {\em frontier} that contain larger numbers of attack-free sensors implicitly acts as a heuristic that shrinks the space to be explored. {This allows our method to significantly outperform competing methods in practice, even for larger problems.}} What's more, although IMHOTEP-SMT can detect the true assignments in the simulation,~\citet{shoukry2017secure} doesn't provide the optimality guarantee.

%% Regarding MIQCP, branch and bound method is used to solve a sequence of 0-1 relaxation problems, each one of which is a Quadratically Constrained Program (QCP). As the number of sensors $p$ increases,  more branches are generated and a larger number of QCP needs to be solved. Hence, MIQCP scales worse when $p$ increases. {In our method, only part of constraints corresponding to a partial assignment are considered when solving the QCP. Therefore, since our algorithm terminates early the search along paths containing invalid nodes, QCPs with large numbers of constraints are typically avoided. Even when the set {\em frontier} contains nodes with higher levels (as the search progresses), it becomes more likely that the corresponding partial assignments are correct, which increases the scalability of our method.}
\section{Conclusion}\label{sec:conclusion}
{In this paper, we proposed a new  optimal graph-search algorithm to correctly identify malicious attacks and to securely estimate the states  in large-scale CPS modeled as linear time-invariant systems.} The graph consists of levels, each one containing two nodes capturing a truth assignment of any given sensor, and directed edges connecting adjacent layers only. Then, our algorithm searches the levels of this graph incrementally, favoring directions with attack-free assignments, while actively managing a repository of partially explored paths with early truth assignments, that can be further explored in the future. The combination of search bias and the ability to self-correct allow our graph-search algorithm to reach the optimal assignment fast. We showed that our algorithm is complete and optimal provided that the attack signal is not dominated by process and measurement noise. Moreover, numerical simulations show that our method outperforms existing algorithms both in terms of optimality and execution time.

\appendix
\section{PROOF OF THEOREM~\ref{thm:comp}}\label{appendixa}
  The proof of completeness of Alg.~\ref{alg:sse} can be divided into two steps. First, we show by contradiction, that each node that is associated with a specific partial truth assignment in the set $node.\ccalI$, can be visited at most once. Then, to show  completeness we show that the algorithm can terminate in a finite number of iterations and the queue {\em frontier} is not empty.

  First, assume that a node associated with a given partial assignment, is visited more than once. Since each partial assignment corresponds to a distinct path from the root to that node, we get that the parent of the node that is associated with this partial assignment except for the last value, is also visited more than once. Repeating this argument, we have that the root is visited more than once, which is a contradiction. Thus, every node associated with a given partial assignment is visited at most once, which means that every node is visited no more times than the  maximum number of  truth assignments starting from the root and ending at that node. For a node that is at level $l$, the maximum number for such truth assignments is $2^l$, i.e., the number of truth assignments of the first $l$ sensors.

  Next, since the nodes and edges in the graph are finite, the previous result implies that the algorithm will terminate in a finite number of iterations. Furthermore, the queue {\em frontier} can not be empty when the algorithm terminates. This is because we do not discard any node unless its partial assignment $\ccalI$ is invalid, that is, the corresponding residual violates $\min_{\mathbf{x} \in \rr{R}{n}} \norm{\mathbf{Y}_{\ccalI} - \ccalO_{\ccalI} \mathbf{x}} \leq  \xoverline{w}_{\ccalI} + \sqrt{\epsilon} $. The reason is that we have shown that if all sensors in $\ccalI$ are attack-free, then we have that $\min_{\mathbf{x} \in \rr{R}{n}}\norm{\mathbf{Y}_{\ccalI} - \ccalO_{\ccalI}\mathbf{x}} \leq \xoverline{w}_{\ccalI} + \sqrt{\epsilon} $ when $\epsilon = \epsilon^*$ and~\eqref{equ:attack} is satisfied.  Therefore, violation of this inequality  means that at least one attacked sensor is incorrectly treated as attack-free. Therefore, only the node with an invalid assignment is discarded. If there exists a solution, the algorithm will search the right node eventually. Once it terminates, the output is a feasible solution. Recall that a child node can be added to the queue {\em frontier} or {\em repo} if $\min_{\mathbf{x} \in \rr{R}{n}}\norm{\mathbf{Y}_{\ccalI} - \ccalO_{\ccalI}\mathbf{x}} \leq \xoverline{w}_{\ccalI} + \sqrt{\epsilon} $. Hence, the constraint~\eqref{equ:mina} is met. Moreover, the constraint \eqref{equ:minb} is checked each time before a child node can be added to {\em frontier} or {\em repo} [Alg.~\ref{alg:sse}, line~\ref{alg:line16}]. Thus, each node existing in {\em frontier} must satisfy~\eqref{equ:minb}, completing the proof.

\section{PROOF OF THEOREM~\ref{thm:opt}}\label{appendixb}

   %% %% The proof for optimality consists of two steps. First, we need to prove that the minimal atttack support is unique and indeed the actual attack support, which is true if the sytem is $2\s$-sparviolatse observable by Theorem {\it II.5} in \cite{shoukry2017secure} and Theorem {\it III.2} in \cite{shoukry2016event}. Next, we show that the search algorithm can provide the minimal attack support solution. Here, we just focus on the second step.
  We prove the optimality of Alg.~\ref{alg:sse} using contradiction. Suppose that the algorithm returns a feasible but suboptimal assignment, which is equivalent to say that the suboptimal path $\pi$, is formed before the optimal $\pi^*$. This is because the algorithm terminates as soon as the first feasible path is found. Since the paths $\pi$ and $\pi^*$ are different, they share identical attack assignments up to a level and then differ in their assignments beyond that level (this level could also be the root).
 %% %% the last level $p$. Otherwise, if the last sensor is attacked, since a feasible path should identify all attacked sensors, these two path will assgin 1 to the node on the highest level, thereby these two paths become identical, and the output is optimal; if the last sensor is attacked-free, recall that when expanding a parent node, Algorithm \ref{alg:sse} assgins 0 first then 1 to generate two children. Both child nodes, denoted by $v_{p}^0$ and $v_p^1$ for level $p$ and corresponding Boolean value, will be put in queue {\em frontier} and node $v_p^0$ with 0 will be popped out first from the queue, thus the output is optimal. Particularly, it these two path divert on the first level, they only share the root.  Next,
 Consider the time when the last shared node $v_{l-1}$, at level $l-1$, is selected to be expanded. Its child node, denoted by \vz with superscript 0 for $value = 0$, is on the optimal path $\pi^*$ and \vo is on the suboptimal path $\pi$. This is because when~\eqref{equ:attack} is satisfied, a feasible path can only treat attack-free sensors as attacked by mistake, but the optimal path does not contain such mistakes.

  Let $\tilde{\pi}^*$ denote the  subpath of $\pi^*$ starting with \vz and $\tilde{\pi}$ denote the subpath of $\pi$ starting with \vo. There are four cases as to whether \vz and \vo can be added to the priority queue {\em frontier} or {\em repo}: (i) both \vz and \vo are added to the priority queue {\em frontier}; (ii) both \vz and \vo are added to the reserved queue {\em repo}; (iii) \vo is added to the {\em frontier} queue but \vz is added to the queue {\em repo}; and (iv) \vz is added to the {\em frontier} queue but \vo is added to the queue {\em repo}. We discuss these four cases separately. The idea is that by assumption, the suboptimal path is formed before the optimal path, thus the algorithm has to switch to expand $v_l^1$ eventually. But the following discussion states that even though the algorithm searches the path starting from $v_l^1$, it will switch back to searching the optimal path.
\begin{figure}[t]
  \centering
  \subfigure[]{
    \label{fig:case1}
    \includegraphics[width=0.45\linewidth]{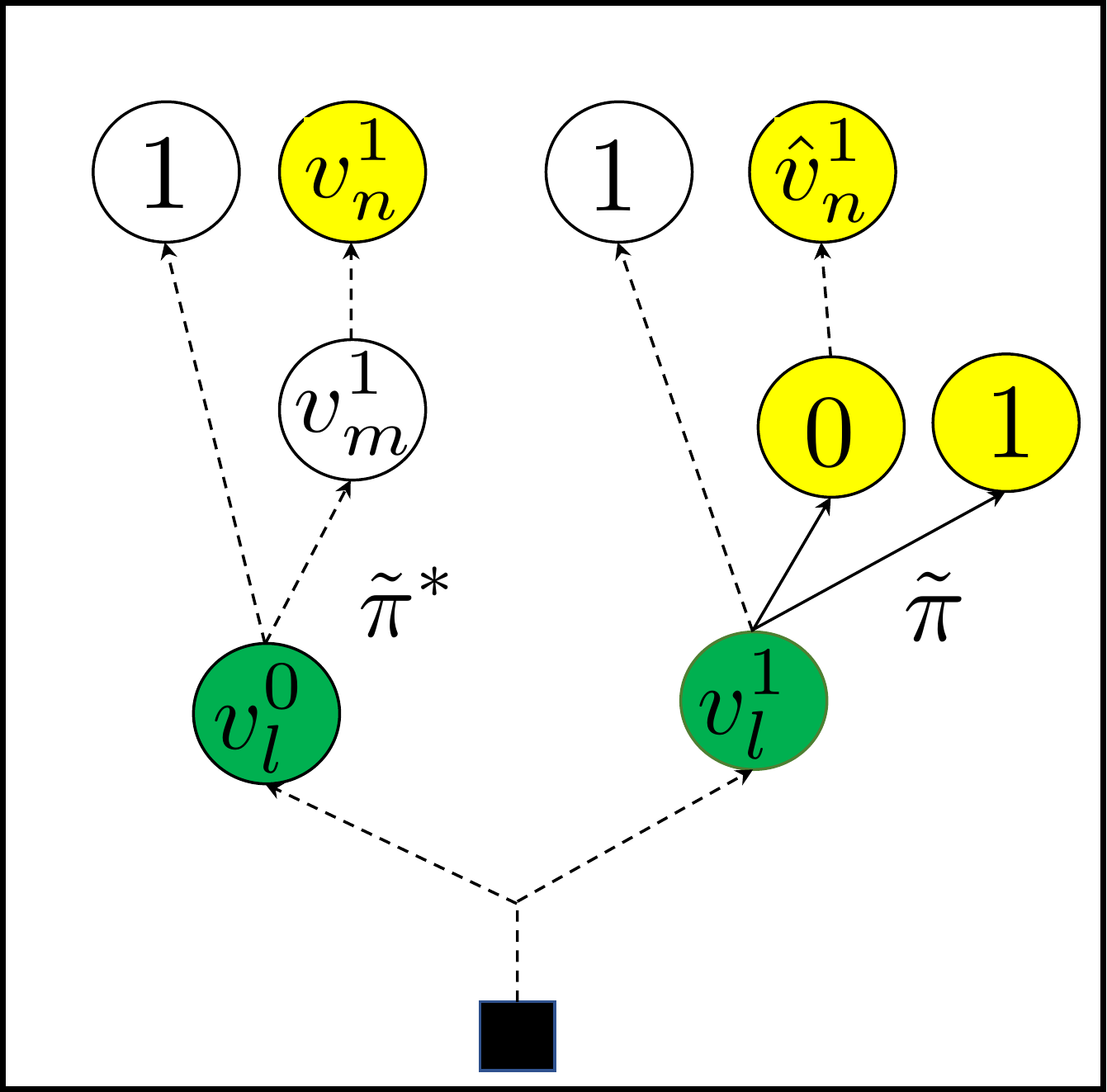}}
   \subfigure[]{
    \label{fig:case3}
    \includegraphics[width=0.45\linewidth]{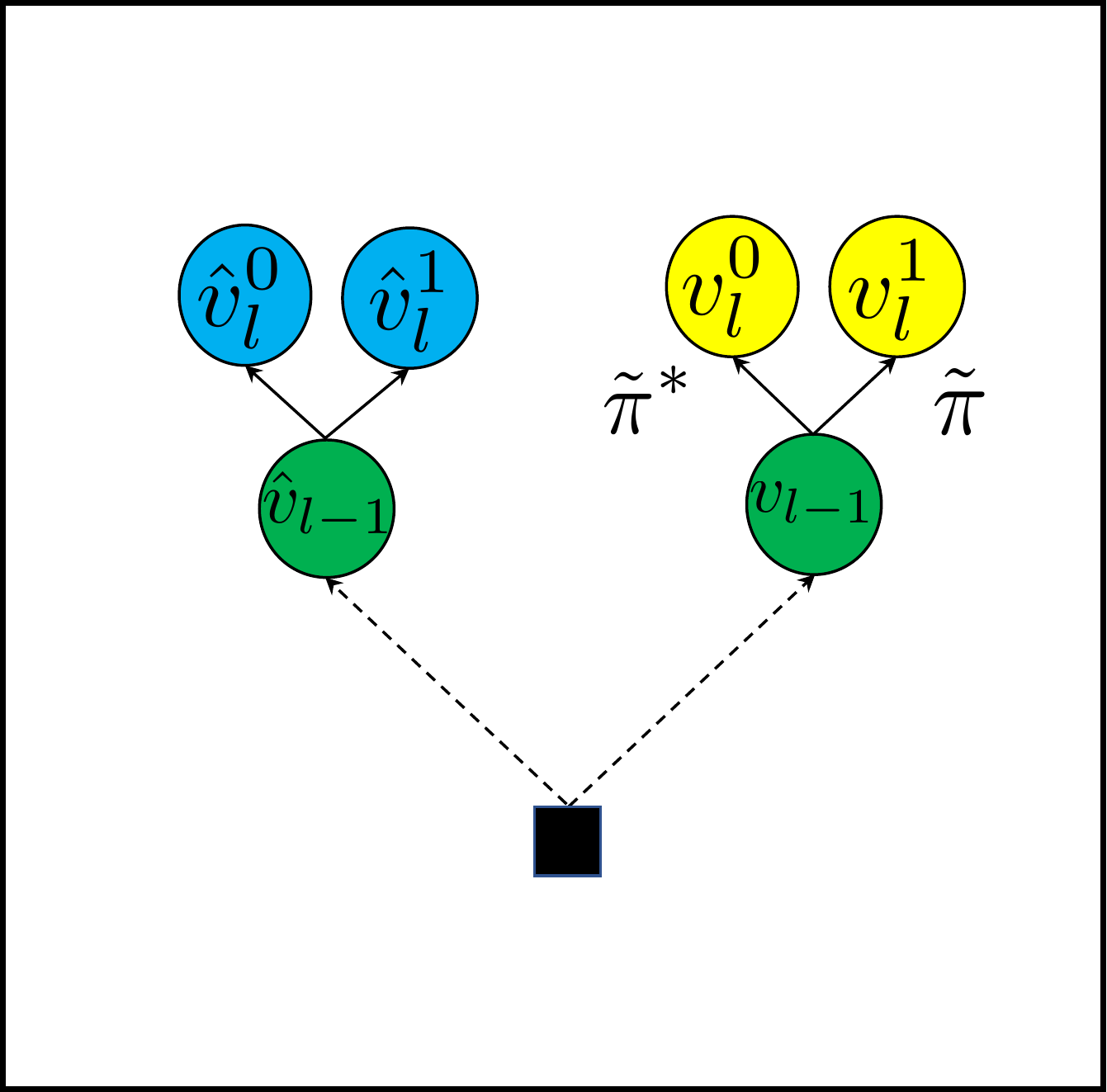}}
   \subfigure[]{
    \label{fig:case4}
    \includegraphics[width=0.45\linewidth]{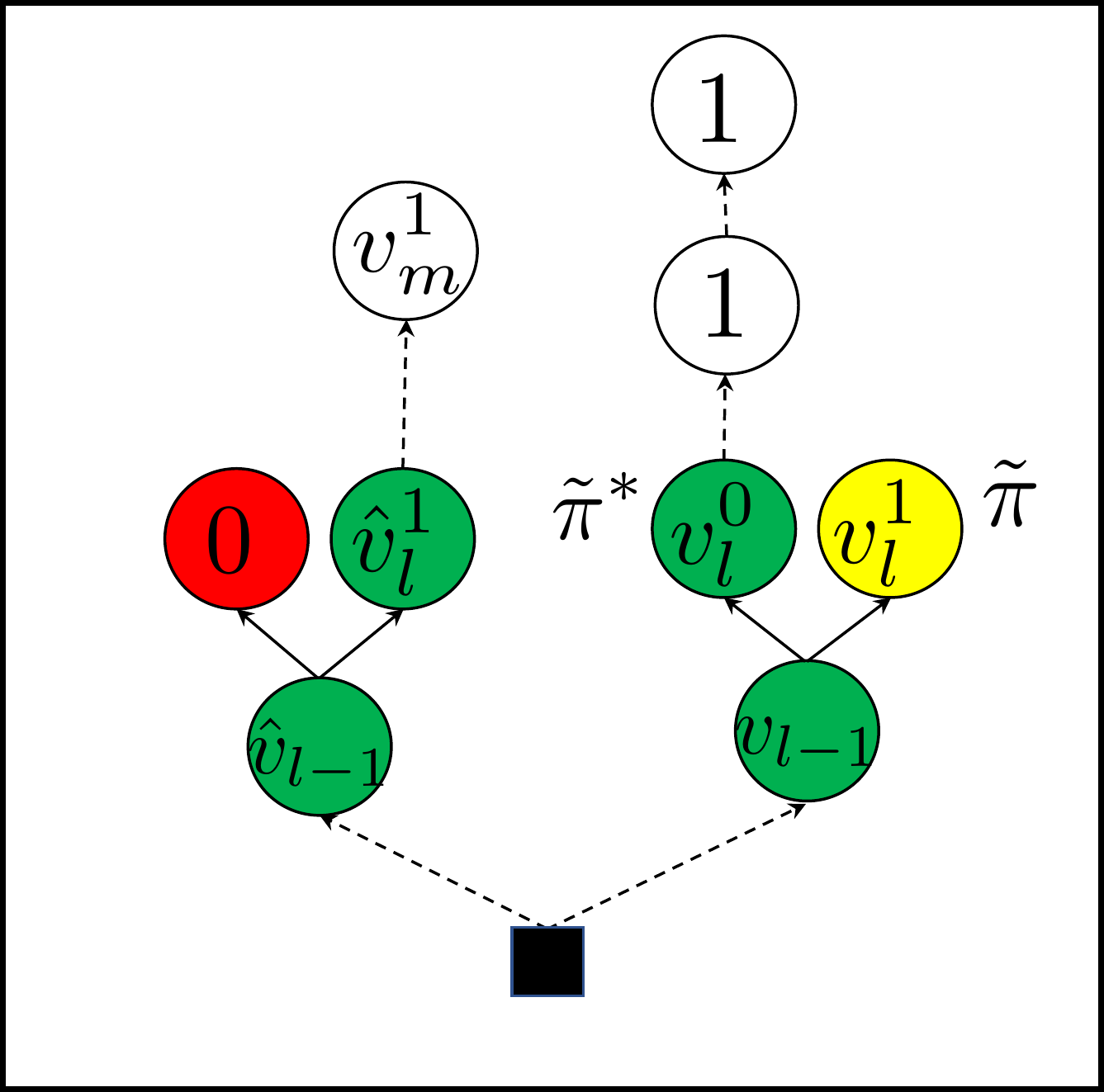}}
  \caption{(a) Case (i), the dashed arrow represents that the tail node is on the path starting from the head node, with intermediate nodes omitted. $v_n^1$ enters {\em repo} because another node with $value = 1$ is in {\em frontier} or {\em explored}, shown as a node with number 1 inside on its left. (b) Case (iii), nodes $\hat{v}_l^0$ and $\hat{v}_l^1$ are in {\em frontier} or {\em explored};  nodes ${v}_l^0$ and ${v}_l^1$ are in {\em repo}; (c) Case (iv), all nodes between $\hat{v}_l^1$ and $v_m^1$ are assigned 1. The level of the first node with $value=1$ on subpath $\tilde{\pi}^*$ can be higher, equal or lower than the level of $v_m^1$.}
\end{figure}

  (i) Both \vz and \vo are added to the priority queue {\em frontier}. Since both \vz and \vo are added to the set {\em frontier}, there does not exist a third node in {\em explored} or {\em frontier} with the same $value$ and $level$ as \vz or $v_l^1$, see Fig.~\ref{appendixb}.\ref{fig:case1}. Therefore, no node with $level = l+1$ is in {\em frontier} or {\em explored}.  Moreover, \vz is selected before \vo since it has higher priority. Expanding $v_l^0$, our algorithm generates two children $v_{l+1}^0$ and $v_{l+1}^1$. If the $(l+1)$-th sensor is attack-free, then both these children nodes will enter {\em frontier}. The algorithm will continue expanding nodes on the path starting from $v_l^0$. When the search algorithm reaches the first node with $value = 1$ on $\tilde{\pi}^*$, denoted by $v_m^1$, then either $v_m^1$ will enter {\em frontier} and  Alg.~\ref{alg:sse} will continue expanding it, or $v_m^1$ will enter {\em repo} and the algorithm will expand $v_l^1$. Since, $v_l^0$ has been expanded, the children of $v_l^1$ will enter {\em repo}. Moreover, $v_m^1$ will be selected from {\em repo} before the children of $v_l^1$ since it has higher priority. Therefore, the algorithm continues searching on the path starting from $v_l^0$. {Since we assume that the assignment corresponding to $\pi$ is the output, the algorithm should switch to searching the children of $v_l^1$ eventually.} In order to expand the children of $v_l^1$, except for $v_m^1$, there must exist another node on $\tilde{\pi}^*$ with $value$ 1, denoted by $v_n^1$, and it should be put in {\em repo}.  This is because if $v_n^1$ enters {\em frontier}, the search will still advance on $\tilde{\pi}^*$. But if $v_n^1$ enters {\em repo}, the children of $v_l^1$ will exit {\em repo} before $v_n^1$, since they contain fewer number of attacked sensors. Meanwhile, the set {\em explored} is emptied. However, the counterpart $\hat{v}_n^1$ of $v_n^1$ on $\tilde{\pi}$ will also enter {\em repo}, since the algorithm will take identical steps searching $\tilde{\pi}$ as when searching $\tilde{\pi}^*$. But $v_n^1$ has higher priority than $\hat{v}_n^1$, therefore, the algorithm switches to searching $\tilde{\pi}^*$ again. Following this logic, we conclude that the search on $\tilde{\pi}$ never surpasses that on $\tilde{\pi}^*$. Intuitively, for those nodes that are assigned 1 on $\tilde{\pi}^*$, the nodes on $\tilde{\pi}$ will be assigned 1 as well, since all feasible paths should detect all attacked sensors correctly. Furthermore, the zero assignments on $\tilde{\pi}$ are a subset of those on $\tilde{\pi}^*$. The nodes on $\tilde{\pi}^*$ always have the advantage over their counterparts on $\tilde{\pi}$ with the same $value$ and $level$ in terms of priority since the optimal path has the least number of attacked sensors. In addition, \vz is expanded before $v_l^1$. Therefore, the last node on the optimal path will exit {\em frontier} before the last node on the suboptimal path. Thus, it is impossible that the suboptimal path is fully formed before the optimal one.

(ii) Both \vz and \vo are added to the reserved queue {\em repo}. If \vz is selected from {\em repo} before $v_l^1$, {\em frontier} only contains \vz and {\em explored}  is empty. Thus, the algorithm starts searching nodes on the path starting from $v_l^0$, as if $v_l^0$ is the new root. As we have discussed in (i), the algorithm will continue searching on the path starting from $v_l^0$. {In order to expand $v_l^1$, since $v_l^1$ is on the suboptimal path that corresponds to the output of the algorithm}, there must exist a second node $v_n^1$ with $value=1$ which  enters {\em repo}. Then $v_l^1$ can be selected from  {\em repo}. But even in this case, when $v_l^1$ is selected for expansion,  the counterpart $\hat{v}_n^1$ of $v_n^1$ on the subpath starting from $v_l^1$ still enters {\em repo}. Same as in (i), the search on $\tilde{\pi}$ never surpasses that on $\tilde{\pi}^*$

   (iii) \vo is added to the {\em frontier} queue but \vz is added to the queue {\em repo}. The graphical illustration of this case is shown in Fig.~\ref{appendixb}.\ref{fig:case3}. If \vz is added to the queue {\em repo}, then a node $\hat{v}_l^0$ with the same $value$ and $level$ as \vz  is already in the queue {\em frontier} or the set {\em explored}. Therefore,  the parent $\hat{v}_{l-1}$ of node $\hat{v}_l^0$, with same level as $v_{l-1}$, is expanded before $v_{l-1}$. Note that only when $\hat{v}_{l-1}$ has higher or equal priority than $v_{l-1}$, can $\hat{v}_{l-1}$ be expanded before $v_{l-1}$.  The reason is that if $v_{l-1}$ has higher priority than $\hat{v}_{l-1}$, then the only way that $\hat{v}_{l-1}$ can be expanded before $v_{l-1}$ is if $v_{l-1}$ or one of its predecessors on $\pi^*$ are in the queue {\em repo} and remain there until they are selected. When this happens, the set {\em explored} is emptied and {\em frontier} contains only $v_{l-1}$ itself or its predecessor. Therefore, there is no need to check whether an equivalent node is in {\em frontier} or {\em explored}, and $v_l^0$ will not enter {\em repo}. Thus, $v_{l-1}$ can not have higher priority than $\hat{v}_{l-1}$. Since $\hat{v}_{l-1}$ is expanded before $v_{l-1}$, $\hat{v}_l^1$ enters {\em frontier} or {\em repo} (because an equivalent node is in {\em frontier} or {\em explored}) before $v_{l-1}$ is expanded. This means that $v_l^1$ will be put in {\em repo}, a contradiction. Hence,
case (iii) is impossible.

  (iv) \vz is added to the {\em frontier} queue but \vo is added to the queue {\em repo}. Like case (iii), since $v_l^1$ is put in {\em repo} when expanding $v_{l-1}$, there exists  a node $\hat{v}_{l-1}$ that is expanded before $v_{l-1}$ and $\hat{v}_l^1$ is in {\em frontier} or {\em explored}. Since $v_l^0$ can enter frontier, we have that $\hat{v}_l^0$ is invalid, see Fig.~\ref{appendixb}.\ref{fig:case4}. $\hat{v}_l^0$ being invalid means that the observability matrix corresponding to the path leading to $\hat{v}_{l-1}$ has full rank or sensor $l$ is under attack. However, the true attack assignment to sensor $l$ is 0, thus the first case holds. Then, $\hat{v}_l^0$ makes the system of linear equations overdetermined and inconsistent, {meaning an attacked sensor is mistreated as attack-free on the path leading to $\hat{v}_l^0$.} Thus, for the following nodes on the path starting from $\hat{v}_l^1$ that shares the same parent $\hat{v}_{l-1}$ with $\hat{v}_l^0$, $value$ can only be set to 1 when those nodes are generated. Assume that when Alg.~\ref{alg:sse} expands $v_{l-1}$, the highest node that is in {\em frontier} or {\em repo} on the path starting from $\hat{v}_l^1$ is $v_m^1$. {Since the algorithm switches to search $v_{l-1}$, we have that $v_{l-1}$ has fewer attacked sensors than $v_m^1$. } After the algorithm switches back to search $v_{l-1}$, we analyze what will happen next by the relationship between the level of the first node with $value=1$ on $\tilde{\pi}^*$ and the level of $v_m^1$. (a) If the level of the first node with $value$ 1 on $\tilde{\pi}^*$ is less than or equal to the level of $v_m^1$, this node will enter {\em repo} and exit {\em repo} before $v_l^1$. (b) If  the first  node on $\tilde{\pi}^*$ with $value$ 1 is at a higher level than $v_m^1$, then the search will advance on  $\tilde{\pi}^*$ until reaching the second node with $value$ 1. But this node will enter {\em frontier} instead of {\em repo}, hence, the search will continue on $\tilde{\pi}^*$. The search on $\tilde{\pi}$ never surpasses that on $\tilde{\pi}^*$.

We conclude that, except case (iii) which does not occur, in the other three cases the algorithm always switches back to expanding the optimal path. This means that the optimal path is formed first, which contradicts our assumption that the suboptimal path is formed before the optimal one, completing the proof.

\bibliographystyle{agsm}
\bibliography{sse}

\end{document}